\documentclass{article} 
\usepackage[left = 2.5cm, right = 2.5cm, bottom = 3cm, top = 2.5cm]{geometry}

\usepackage{amsmath,amssymb,amsthm,bm,bbm,dsfont}

\DeclareMathOperator{\expect}{E}
\DeclareMathOperator{\var}{var}

\usepackage{booktabs}
\usepackage{multirow}
\usepackage{rotating}
\usepackage{graphicx}
\usepackage{subcaption}
\usepackage{tabularx}
\usepackage{dcolumn}

\usepackage{pdfpages}

\usepackage{natbib}
\usepackage{hyperref}
\usepackage{xr}
\usepackage{autonum}
\usepackage{siunitx}

\usepackage{enumitem}

\usepackage[framemethod=TikZ]{mdframed}

\newtheorem{theorem}{Theorem}[section]
\newtheorem{conjecture}{Conjecture}[section]

\theoremstyle{definition}

\surroundwithmdframed[
topline=false,
rightline=false,
bottomline=false,
leftmargin=\parindent,
skipabove=\medskipamount,
skipbelow=\medskipamount,
linecolor=grey!30,
innerlinewidth=1pt
]{proof, definition, remark, theorem, lemma, proposition, corollary}

\makeatletter
\newcommand{\vnorm}[1]{%
	\ensuremath{%
		\if@display
		\left\| #1 \right\|
		\else
		\| #1 \|
		\fi
	}%
}
\makeatother
\newcommand{\mnorm}[1]{{\left\vert\kern-0.25ex\left\vert\kern-0.25ex\left\vert #1 
		\right\vert\kern-0.25ex\right\vert\kern-0.25ex\right\vert}}
\newcommand{\mnorms}[1]{{\vert\kern-0.25ex\vert\kern-0.25ex\vert #1 
		\vert\kern-0.25ex\vert\kern-0.25ex\vert}}
\newcommand*{\bb}{\boldsymbol}

\newcommand{\ttextrm}[1]{\textrm{\tiny #1}}

\newcommand{\prox}[2]{\ensuremath{\textrm{prox}_{#1}\left(#2\right)}}
\newcommand{\betady}{\ensuremath{\hat{\bbeta}^{\ttextrm{DY}}}}
\newcommand{\etady}{\ensuremath{\hat{\eta}^{\ttextrm{DY}}}}
\newcommand{\setady}{\ensuremath{\tilde{\eta}^{\ttextrm{DY}}}}
\newcommand{\thetady}{\ensuremath{\hat{\btheta}^{\ttextrm{DY}}}}
\newcommand{\betar}{\ensuremath{\hat{\bbeta}^{\ttextrm{R}}}}
\newcommand{\betam}{\ensuremath{\hat{\bbeta}^{\ttextrm{ML}}}}

\def\bW {\bb{W}}
\def\bX{\bb{X}}
\def\bx{\bb{x}}

\def\b0{\bb{0}}

\def\bv{\bb{v}}
\def\by{\bb{y}}
\def\bA{\bb{A}}

\def\bSigma{\bb{\Sigma}}

\def\bbeta{\bb{\beta}}

\def\btheta{\bb{\theta}}

\def\mus{{\mu}_{*}}
\def\bs{{b}_{*}}
\def\sigmas{{\sigma}_{*}}

\usepackage{xcolor}

\usepackage{authblk}
\usepackage{orcidlink}

\author[1]{Philipp Sterzinger~\orcidlink{0009-0007-7348-5810}\thanks{p.sterzinger@lse.ac.uk}}
\author[2]{Ioannis Kosmidis~\orcidlink{0000-0003-1556-0302}\thanks{ioannis.kosmidis@warwick.ac.uk}}

\affil[2]{Department of Statistics, London School of Economics \authorcr London, WC2A 2AE, UK}
\affil[2]{Department of Statistics, University of Warwick \authorcr Coventry, CV4 7AL, UK}

\title{Diaconis-Ylvisaker prior penalized likelihood for $p/n \to \kappa \in (0,1)$ logistic regression}

\begin{document}
\maketitle 

\begin{abstract}
  We characterise the behavior of the maximum Diaconis--Ylvisaker
  prior penalized likelihood estimator in high-dimensional logistic
  regression, where the number of covariates is a fraction
  $\kappa \in (0,1)$ of the number of observations $n$, as
  $n \to \infty$. We construct a rescaled estimator with zero
  asymptotic aggregate bias and define adjusted $Z$-statistics and
  rescaled penalized likelihood ratio statistics that exhibit the
  typical null asymptotic distributions, when the covariates are
  independent multivariate normal with an arbitrary covariance matrix
  and the linear predictor has asymptotic variance $\gamma^2$. While
  the maximum likelihood estimate asymptotically exists only for a
  narrow range of $(\kappa, \gamma)$ values, the maximum
  Diaconis--Ylvisaker prior penalized likelihood estimate always exists
  and can be computed directly using standard maximum likelihood
  routines. Thus, our asymptotic results extend to $(\kappa, \gamma)$
  values where the maximum likelihood framework breaks down, with no
  additional implementation or computational cost. We study the
  estimator's shrinkage properties, compare the proposed estimation
  and inference procedures with alternatives that also accommodate
  proportional asymptotics, and formulate a conjecture --- supported
  by strong empirical evidence --- that extends our results when the
  model includes an intercept parameter. Finally, we propose
  estimation methods for all unknown constants involved in our
  procedures and demonstrate the theoretical advances through
  extensive simulation studies and the analysis of digit
  recognition data.
\end{abstract}

\section{Introduction} 

\subsection{Logistic regression} 

Logistic regression is arguably one of the most widely used models in
statistical practice to associate binary responses with a sequence of
covariates, for either inference on covariate effects or prediction. A
logistic regression model assumes that, conditionally on covariate
vectors $\bx_1,\ldots,\bx_n$, $\bx_i \in \Re^p$, the responses
$y_1,\ldots,y_n$, $y_i \in \{0,1\}$ are realisations of independent
Bernoulli random variables such that
\begin{align}
  \label{eq:logistic_y}
  \Pr(y_i = 1 \mid \bx_i) = \zeta'\left(\bx_i^\top \bbeta_0\right) \quad (i = 1,\ldots,n) \,, 
\end{align}
for some unknown vector of regression coefficients
$\bbeta_0\in \Re^p$, and where $\zeta'(x) = 1 / ({1+e^{-x}})$ is the
derivative of $\zeta(x) = \log(1+e^x)$. The log-likelihood of the
logistic regression model is
$\ell(\bbeta; \by, \bX) = \sum_{i=1}^{n} \left\{y_i \bx_i^\top \bbeta
  - \zeta\left(\bx_i^\top \bbeta\right)\right\}$, where
$\by = (y_1,\ldots,y_n)^\top$ and $\bX$ is the $n \times p$ matrix
with $j$th row $\bx_i$. The maximum likelihood (ML) estimator $\betam$
is the maximiser of $\ell(\bbeta; \by, \bX)$. Despite the widespread
use of ML for estimating the signal $\bbeta_0$, the existence of the
ML estimate is not guaranteed for all configurations of responses and
covariates. Specifically, \citet{albert+anderson:1984} show that the
ML estimate does not exist --- colloquially, at least one of its
components takes on an infinite value --- if and only if there is data
separation. Data separation occurs when there exists a separating
hyperplane in $\Re^p$ that perfectly discriminates among the
covariates $\bx_1,\ldots,\bx_n$ according to the values of the
corresponding outcomes $y_1,\ldots,y_n$, except, perhaps, from
observations falling onto that separating
hyperplane. \citet{candes+sur:2020} show that in high-dimensional
logistic regression with $p/n \to \kappa \in (0,1)$,
$\bx_i \sim \mathrm{N}(\bb{0}_p, \bSigma)$ and signal strength
$\var(\bx_i^\top \bbeta_0) \to \gamma^2 \in (0, \infty)$ with,
potentially, non-zero intercept parameter, the ML estimate exists with
probability one only in a narrow region under a phase transition curve
in the $(\kappa, \gamma)$-plane.


\subsection{Maximum penalized likelihood estimation} 

There are several maximum penalized likelihood (MPL) approaches that
can remedy the occurrence of infinite ML estimates in logistic
regression. These methods maximise a penalized log-likelihood
$\ell(\bbeta; \by, \bX) + \log p(\bbeta)$ where the penalty
$\log p(\bbeta)$ is chosen to ensure that the MPL estimate
exists. Simple separable functions of the form
$\log p(\bbeta) = \sum_{j=1}^p f(\bbeta_j)$ are popular penalty
choices; see, for example \citet[Section~3.2]{hastie+etal:2015} for
the LASSO and \citet{cessie+houwelingen:1992} for logistic ridge
regression. Notable non-separable penalties that depend on $\bX$ and
deliver finite estimates are the logarithm of the Jeffreys' invariant
prior \citep[see][Corollary~1]{kosmidis+firth:2021} and the logarithm
of the conjugate prior of \citet{diaconis+ylvisaker:1979}, henceforth
called DY prior \citep[see][Theorem~1]{rigon+aliverti:2023}, which is
the subject of investigation in the current work. The log of the DY
prior for the logistic regression model is given by
\begin{equation}
  \label{eq:prior}
  \log p(\bbeta;\bX) = \frac{1 - \alpha }{\alpha} \sum_{i=1}^{n}\left\{ \zeta'\left(\bx_i^\top  \bbeta_P\right) \bx_i^\top \bbeta  - \zeta\left(\bx_i^\top \bbeta\right) \right\} + C \,, 
\end{equation} 
where $C$ is a normalising constant that does not depend on $\bbeta$,
$\bbeta_{P} \in \Re^p$ is the mode of the DY prior, and
$\alpha \in (0, 1]$ controls the variability of the prior
distribution. \citet{diaconis+ylvisaker:1979} provide further details
and discussion from a Bayesian viewpoint.

In the classical fixed-$p$ setting, the asymptotic properties of the
maximum Jeffreys’ prior penalized likelihood (MJPL) and maximum
Diaconis-Ylvisaker prior penalized likelihood (MDYPL) estimators are
well understood. Under mild regularity conditions, the bias of the
MJPL estimator decreases faster with the amount of information about
the model parameters than the bias of the ML estimator, and both the
MJPL and MDYPL estimators preserve the asymptotic normality and
Cram\'{e}r-Rao efficiency that is typically expected of the ML
estimator. Furthermore, the MJPL estimator is equivariant only under
linear transformations of the parameters, while the MDYPL estimator,
inherits all equivariance properties the ML estimator has
\citep[see][for details]{zehna:1966}. \citet{kosmidis+firth:2021} and
\citet{rigon+aliverti:2023} provide details about MJPL and MDYPL,
respectively.

\subsection{Crossing the phase transition}
\label{sec:phase_transition}

\begin{figure}[t]
  \centering
  \includegraphics[width = \linewidth]{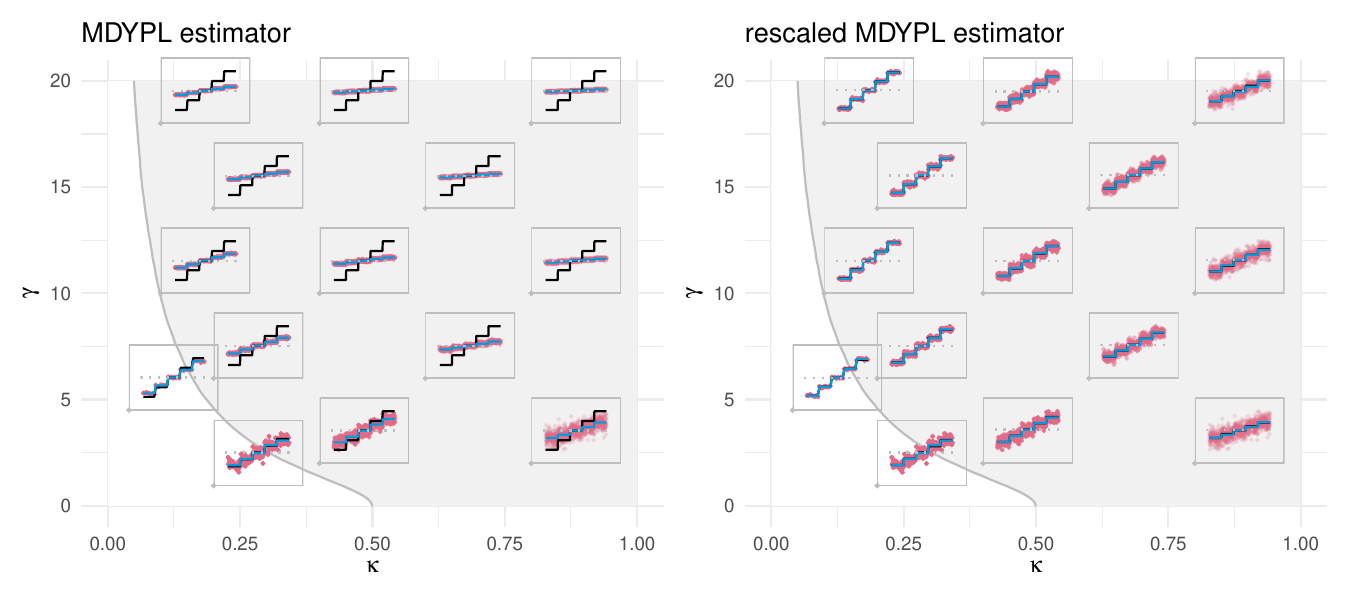}
  \caption{MDYPL (left) and rescaled MDYPL (right) estimates for
    various configurations of $(\kappa,\gamma)$ and
    $\alpha = 1 / (1+\kappa)$ in the simulation setting of
    Section~\ref{sec:phase_transition}. The white and grey area
    indicate where the ML estimate does or does not exist
    asymptotically, respectively. Red points show the average
    coefficient estimates over $10$ independent replications per
    $(\kappa, \gamma)$ setting. The cyan segments show the sample mean
    of the estimates for each value of the truth, and the black
    segments are the truth.}
  \label{fig:ra_scaled}
\end{figure}   

The guaranteed existence of MJPL and MDYPL estimates makes these
estimators natural candidates for high-dimensional regression
problems. Recent empirical results show that the MJPL estimator
performs well without adjustments in high-dimensional settings where
ML estimates exist and should be rescaled whenever the ML estimator
does not exist \citep[see, for example,][]{kosmidis+zietkiewicz:2025},
and \citet{rigon+aliverti:2023} illustrate that the MDYPL estimator
with $\bbeta_{P} = \b0_p$ and $\alpha = 1 / (1 + \kappa)$ in
expression~\eqref{eq:prior} has good empirical performance in terms of
bias and mean squared error (MSE) for specific fixed-$p$ regimes, and
in a setting with low signal strength ($\gamma = \sqrt{0.9}$) and
moderate dimension ($\kappa = 0.2$).

Apart from these empirical results, little is known about the
behaviour of MJPL and MDYPL in high dimensions with
$p/n \to \kappa \in (0, 1)$. For a more thorough investigation of the
performance of the MDYPL estimator, we extend the computer experiment
in \citet[Section~4.3]{rigon+aliverti:2023} to a range of
$(\kappa, \gamma)$ values.  We consider $n = 1000$ independent
covariate vectors with $\bx_i\sim \mathrm{N}(\b0_p, n^{-1}\bb{I}_p)$,
where $\bb{I}_p$ is the $p \times p$ identity matrix, set $\bbeta_0$
to $p/5$ replications of the vector $(-3,-3/2,0,3/2,3)^\top$ with
$p = n \kappa$, scaled to satisfy
$\var(\bx_i^\top \bbeta_0) = \gamma^2$, for given $\kappa$ and
$\gamma$.

Figure~\ref{fig:ra_scaled} shows MDYPL estimates, averaged over ten
independent replications at various points in the $(\kappa,\gamma)$
plane, including $(0.2, \sqrt{0.9})$.  The MDYPL estimate exists for
$(\kappa, \gamma)$ values beyond the phase transition curve of
\citet{candes+sur:2020}. In addition, the good empirical
performance of the MDYPL estimator for $\kappa = 0.2$ and
$\gamma = \sqrt{0.9}$ is not uniform across the
$(\kappa, \gamma)$-plane, and that performance degrades substantially
beyond the phase transition curve. Indeed, the MDYPL estimates are
practically useless for recovering the signal $\bbeta_0$ for the
majority of the $(\kappa, \gamma)$-plane, and demonstrate severe
shrinkage towards zero as $\kappa$ and $\gamma$ increase. The right
panel of Figure~\ref{fig:ra_scaled} shows that rescaling the MDYPL
estimates according to the developments in the current work reliably
recovers the signal across the points on the $(\kappa, \gamma)$-plane
we consider.

\subsection{Approximate message passing for logistic regression}



\citet{donoho+etal:2009} introduced approximate message passing (AMP)
as a class of iterative algorithms that take as input a random matrix
$\bA \in \Re^{n \times p}$ and output some iterates of interest
$x^t \in \Re^p, z^t \in \Re^p$. These algorithms have the remarkable
property that their iterates have a tractable asymptotic limiting
behaviour in the high-dimensional limit where both $p,n \to \infty$,
as first proven in
\citet{bayati+montanari:2011}. \citet{sur+candes:2019} extend the work
of \citet{javanmard+montanari:2013} and \citet{donoho+montanari:2016}
and develop an AMP recursion that
describes the aggregate limiting behaviour of the ML estimator in
high-dimensional logistic regression with no intercept when
$p/n \to \kappa \in (0, 1)$ and the covariates are independent with
$\bx_i \sim \mathrm{N}(\bb{0}_p, n^{-1}\bb{I}_p)$. In particular,
\citet{sur+candes:2019} show that the ML estimator exhibits, in
settings where it exists, bias in an aggregate sense, inflated
standard errors and that the likelihood-ratio test statistic converges
in distribution to a scaled $\chi^2$ random
variable. \citet{zhao+etal:2022} generalise these results to normal
model matrices with arbitrary covariate covariance and provide results
about the asymptotic distribution of the ML estimator's coordinates,
as well as the likelihood ratio test statistic, using a stochastic
representation of the ML estimator. In the same spirit,
\citet{salehi+et+al:2019} establish the aggregate asymptotic behaviour
for maximum penalized likelihood estimators with separable, convex
penalty functions, focusing on the regime where
$p/n \to \kappa \in (0,\infty)$. The results of
\citet{salehi+et+al:2019} are, though, specific to aggregate
asymptotics with no consideration of inference.

\subsection{Our contribution} 

In the current work we: 1) establish the asymptotic aggregate
behaviour of the MDYPL estimator in high-dimensional logistic
regression with independent normal covariates with arbitrary
covariance matrix under the proportional asymptotic regime with
$p /n \to \kappa \in (0, 1)$ and limiting signal strength
$\gamma^2 \in (0, \infty)$; 2) derive adjusted $Z$-statistics and
rescaled DY prior penalized likelihood ratio statistics that recover
inferential performance about the model parameters; 3) provide a
conjecture with strong empirical support about the extension of 1) and
2) for logistic regression with unknown intercept; and 4) present
procedures for the estimation of all unknown constants that describe
the asymptotic behaviour of the MDYPL estimator. We choose MDYPL as
the basis for those developments because it is convenient both for
estimation and inference, essentially being ML with deterministically
adjusted responses prior to fitting.
In contrast to corresponding results about ML in \citet{sur+candes:2019} and
\citet{zhao+etal:2022}, our advances are no longer confined to the
narrow region of $(\kappa, \gamma)$ values below the phase transition
curve of \citet{candes+sur:2020} where the ML estimates asymptotically
exist. We also show that 1)--4) generalise of
the corresponding results about ML estimation, reducing to those for
particular values of the DY prior hyperparameters.

We introduce novel analytical techniques to the AMP-literature that
can find use in other AMP based work. Firstly, while the DY prior
penalty is nonseparable (see the seminal works of
\citealp{berthier+montanari+nguyen:2020} and \citealp{huang:2022} on
nonseparable AMP), the DY prior penalized log-likelihood can be
expressed as a logistic regression log-likelihood with transformed
responses. This enables the study of MDYPL estimation without having
to explicitly accommodate for the nonseparable penalty function, and
paves the way for similar results for other penalty functions which
reduce to a perturbation of the responses. Secondly, we provide a
rigorous treatment of the additional smoothness steps required to
formally transfer the asymptotic results from the AMP literature to
logistic regression, by introducing a Lipschitz-smooth AMP recursion
where the binary responses,
$Y_i = \mathds{1} \left\{\bar{\varepsilon}_i < \zeta'\left(\bx_i^\top
    \bbeta_0 \right)\right\}$ for $\bar{\varepsilon}_i \sim U(0,1)$
--- which are not Lipschitz-smooth in their arguments
$\bar{\varepsilon}_i$ and $\bx_i^\top \bbeta_0$ --- are approximated
by a smooth function. This step has been missing in the analysis of
the ML estimator in high-dimensional logistic regression in
\citet{sur+candes:2019}, as noted in \citet{feng+et+al:2022}. Details
are provided in Section~\ref{subsubsec:smooth_AMP} of the
Supplementary Material document.

We study the MDYPL estimator's shrinkage properties under various
optimality criteria, and compare MDYPL to logistic ridge regression
\citep{salehi+et+al:2019} and corrected-least squares
\citep{lewis+battey:2023}, which are two methods for which the
asymptotic performance of the respective estimators has been examined
in the high-dimensional regime we consider. We accompany our
theoretical findings with extensive simulation studies, and provide a
real-data case study on digit recognition. Overall, we provide a
comprehensive picture of the performance of the MDYPL framework in the
challenging asymptotic regime of proportional asymptotics.

All the methods we develop are implemented in the brglm2 R package
\citep{brglm2}. The proofs of all theoretical results, along with
further empirical and numerical results are provided in the
Supplementary Material document, which is available at
\url{https://github.com/psterzinger/MDYPL}.





\section{Maximum Diaconis-Ylvisaker prior penalized likelihood}
\label{sec:dy_mape}

Ignoring normalising constants, the DY prior penalized log-likelihood
for logistic regression takes the form
\begin{equation}
  \label{eq:posterior}
  \ell^*(\bbeta;\by, \bX) =   \frac{1}{\alpha} \sum_{i=1}^n \left\{ \left(\alpha  y_i + (1 - \alpha) \zeta'\left(\bx_i^\top \bbeta_P\right) \right) \bx_i^\top\bbeta   -  \zeta\left(\bx_i^\top \bbeta\right) \right\} \,. 
\end{equation}
The maximiser of \eqref{eq:posterior} is the MDYPL estimator and is
denoted by $\betady$.
Note that $\ell(\bbeta;\by^*,\bX) / \alpha = \ell^*(\bbeta; \by, \bX)$, where
$\by^*= (y_1^*, \ldots, y_n^*)^\top$, with
$y_i^{*} = \alpha y_i + (1-\alpha)\zeta'(\bx_i^\top \bbeta_P)$. Hence,
$\betady$ is the maximiser of a logistic regression log-likelihood
with pseudo-responses $\by^*$, and MDYPL estimates can be computed
using standard ML routines for logistic regression. In contrast to the
ML estimate, the MDYPL estimate is unique and exists for all data
configurations $\{\by,\bX\}$, with full rank $\bX$; see, for
example, \citet[Theorem~1]{rigon+aliverti:2023}.

The MDYPL estimator converges to the ML estimator for $\alpha \to 1$,
and to the prior mode $\bbeta_P$ for $\alpha \to 0$. For this reason,
we refer to $\alpha$ as shrinkage parameter. The ML estimator of a
logistic regression model is biased away from $\bb{0}_p$
\citep[see][Section~8]{cordeiro+mccullagh:1991} in the fixed-$p$
asymptotic regime, which appears to also be the case in high
dimensions according to recent empirical evidence (see
\citealt{sur+candes:2019} and \citealt{zhao+etal:2022}). Hence, it is
natural to shrink the MDYPL estimator away from the ML estimator
towards zero by setting $\bbeta_{P} = \bb{0}_p$.  In what follows, all
results we derive hold for $\bbeta_{P} = \bb{0}_p$ and any
$\alpha \in (0, 1)$.


\section{Asymptotic behaviour of the MDYPL estimator}
\label{sec:main}

\subsection{Preamble}

All theoretical results are developed under the same assumptions
behind the advances for ML estimation in \citet{sur+candes:2019} and
\citet{zhao+etal:2022}: (i) \textit{Proportional asymptotics}:
We operate in an asymptotic regime where the number of covariates $p$
can be comparable to the number of observations $n$, in the sense that
$p / n \to \kappa \in (0,1)$. (ii) \textit{Normal covariates}: The
rows of the model matrix $\bX$ are independent and identical
realisations of multivariate normal distributions with mean zero. The
simulation results in Section~\ref{sec:case-study} and the empirical
evidence in Section~\ref{sec:bt-mar} of the Supplementary Material
document on the estimation and inference from Bradley-Terry models
suggest that this assumption can be relaxed to sub-Gaussian
covariates. (iii) \textit{Fixed asymptotic signal strength}: The
variance of the linear predictors at the signal $\bbeta_{0}$ satisfies
$\var(\bx_i^\top \bbeta_{0}) \to \gamma^2$ for
$\gamma \in (0, \infty)$. This ensures that $\bx_i^\top \bbeta_0$ does
not tend to $0$ or diverge to $\infty$ and that the problem is truly
high-dimensional ($\gamma \neq 0$) and non-trivial, in that the
likelihood does not converge to $0$ or $1$.

\subsection{Aggregate behaviour}
\label{subsec:aggregate}

Central to the developments is the asymptotic behaviour of the
divergence
$\sum_{j = 1}^p \psi(\betady_i - \mus \bbeta_{0,j}, \bbeta_{0,j}) / p$
as $n \to \infty$, where $\psi: \Re^2 \to \Re$ is a pseudo-Lipschitz
function of order two (see Definition~\ref{def:pl} in the
Supplementary Material document), and $\mus$ is a scalar that we
define in Section~\ref{subsec:aggregate}.  As a starting point, we
derive the almost sure limit of the divergence in the regime of
\citet{sur+candes:2019} where
$\bx_i \sim \mathrm{N}(\b0_p, n^{-1}\bb{I}_p)$ and
$\var\left(\bx_i^\top \bbeta_{0} \right) \to \gamma^2$. In a similar
vein to the work of \citet{zhao+etal:2022}, that result can be
extended to derive the asymptotic behaviour of the divergence with
normal covariates $\bx_i \sim \mathrm{N}(\b0_p,\bSigma)$ and
$\var\left(\bx_i^\top \bbeta_{0}\right) = \bbeta_0^\top \bSigma
\bbeta_0 \to \gamma^2$, under mild regularity conditions on $\bSigma$.

The asymptotic behaviour of $\betady$ is analysed using an
appropriately defined generalised AMP recursion, whose iterates have a
known asymptotic distribution as $p / n \to \kappa \in (0,1)$, with a
stationary point that coincides with the MDYPL estimator. The
asymptotic behaviour of the MDYPL estimator is governed by the
solution $(\mus,\bs,\sigmas)$ to the system of nonlinear equations in
three variables $(\mu,b,\sigma)$
\begin{equation}
  \begin{aligned}
    \label{eq:state_evol}
    \expect \left[2 \zeta'(Z) Z \left\{\frac{1+\alpha}{2} - \zeta' \left(\prox{b \zeta}{Z_{*} + \frac{1+\alpha}{2}b} \right) \right\} \right] &= 0 \\ 
    1 - \kappa - \expect \left[\frac{2 \zeta'(Z)}{1 + b \zeta''(\prox{b\zeta}{Z_{*} + \frac{1+\alpha}{2}b})} \right] &= 0 \\ 
    \sigma^2 -  \frac{b^2}{\kappa^2}  \expect \left[2 \zeta'(Z) \left\{\frac{1+\alpha}{2} - \zeta'\left(\prox{b\zeta}{  Z_{*} + \frac{1+\alpha}{2}b}\right) \right\}^2 \right] &= 0 \,,
  \end{aligned}
\end{equation}
where $Z \sim \mathrm{N}(0,\gamma^2)$,
$Z_{*} = \mu Z + \kappa^{1/2} \sigma G$ for $G \sim \mathrm{N}(0,1)$
which is independent of $Z$, $\zeta''(\cdot)$ is the second derivative
of $\zeta(\cdot)$, and
$\prox{\bs \zeta}{x}=\arg \min_{u} \left\{ \bs \zeta(u) +
  \frac{1}{2}(x-u)^2 \right\}$ denotes the proximal operator.  The
system of equations in \eqref{eq:state_evol} defines the stationary
points to the state evolution of the underlying AMP recursion, which
can also be derived as the first-order optimality conditions of a
related auxiliary optimisation problem to MPL estimation (see, for
example, \citealt[Appendix~C]{salehi+et+al:2019}). The state evolution
equations require knowledge of the unknown asymptotic signal strength
$\gamma$ and the dimensionality constant $\kappa$, which, in practice,
must be estimated from the data. We discuss the estimation of these
unknown constants in Section~\ref{sec:unknown_constants}.

By properties of the proximal operator and $\zeta(\cdot)$,
manipulations similar to those used in Lemma~\ref{lemma:eq_state_ev}
of the Supplementary Material document show that the system of
equations in \eqref{eq:state_evol} coincides with
\citet[equation~(5)]{sur+candes:2019}, for $\alpha = 1$. As a
result, Theorem~\ref{thm:aggregate_asymptotics} below encompasses the
asymptotic behaviour of the ML estimator, in settings where the ML estimate
exists. Denote by $\bb{J}(\mu,b,\sigma)$ the Jacobian matrix of the
LHS of \eqref{eq:state_evol} with respect to $(\mu,b,\sigma)$.
\begin{theorem}
  \label{thm:aggregate_asymptotics}
  Consider the logistic regression model~(\ref{eq:logistic_y}) with
  independent covariates
  $\bx_i \sim \mathrm{N}(\bb{0}_p,n^{-1}\bb{I}_p)$, and the empirical
  distribution of the elements of $\bbeta_0$ converging weakly to
  $\bar{\beta} \sim \pi_{\bar{\beta}}$ with
  $ \sum_{j = 1}^p \bbeta_{0,j}^2 / p \to \expect(\bar{\beta}^2) <
  \infty$. Assume that
  $(\alpha,\kappa,\gamma)$ are such that \eqref{eq:state_evol} admits
  a solution $(\mus,\bs,\sigmas)$ such that
  $\bb{J}(\mus,\bs ,\sigmas)$ is nonsingular. Then, for any function $\psi: \Re^2 \to \Re$ that is
  pseudo-Lipschitz of order 2,
  \begin{equation}
    \label{eq:thm_eq}
    \frac{1}{p} \sum_{j=1}^{p} \psi(\betady_j - \mus \bbeta_{0,j}, \bbeta_{0,j}) \overset{\textrm{a.s.}}{\longrightarrow} \expect \left[\psi(\sigmas G, \bar{\beta})\right], \quad \textrm{as } n \to \infty \,, 
  \end{equation}
  where $G \sim \mathrm{N}(0,1)$ is independent of $\bar{\beta}$.   
\end{theorem}

For appropriate choices of $\psi(t, u)$, one gets the typical
asymptotic results for centred AMP recursions: a) $p^{-1} \sum_{j=1}^p \left(\betady_j - \mus \bbeta_{0,j}\right) \overset{\textrm{a.s.}}{\longrightarrow} 0$, for $\psi(t,u) = t$; b) $p^{-1}\sum_{j=1}^p \left(\betady_j - \mus \bbeta_{0,j}\right)^2 \overset{\textrm{a.s.}}{\longrightarrow} \sigmas^2$, for $\psi(t,u) = t^2$; c) $p^{-1}\vnorm{\betady-\bbeta_0}_2^2 \overset{\textrm{a.s.}}{\longrightarrow} \sigmas^2 + (1-\mus)^2  \gamma^2 / \kappa$, for $\psi(t,u) = (t-(1-\mus)u)^2$; d) $p^{-1}\sum_{j=1}^p (\betady_j - \mus \bbeta_{0,j} ) \bbeta_{0,j} \overset{\textrm{a.s.}}{\longrightarrow} 0$, for $\psi(t,u) = tu$.
The parameter $\mu^*$ in statement a) is termed the aggregate bias
parameter in \citet{sur+candes:2019}, b) and c) characterise
aggregate measures of variance and MSE, and d) shows that the
recentred estimator is asymptotically uncorrelated from the signal
$\bbeta_0$ in an aggregate sense.  Since we are typically interested
in signal recovery, we can also obtain corresponding results for the
rescaled MDYPL estimator $\betady / \mus$.
In particular,
e)
$p^{-1} \sum_{j=1}^p (\betady_j / \mus - \bbeta_{0,j})
\overset{\textrm{a.s.}}{\longrightarrow} 0$, for
$\psi(t,u) = t / \mus$; f)
$p^{-1}\sum_{j=1}^p (\betady_j / \mus - \bbeta_{0,j})^2
\overset{\textrm{a.s.}}{\longrightarrow} \sigmas^2 / \mus^2$, for
$\psi(t,u) = t^2/\mus^2$; g)
$p^{-1}\sum_{j=1}^p (\betady_j / \mus- \bbeta_{0,j} ) \bbeta_{0,j}
\overset{\textrm{a.s.}}{\longrightarrow} 0$, for
$\psi(t,u) = tu / \mus$.  If the empirical distribution of the
components of $\bbeta_0$ converges to a discrete distribution
$\pi_{\bar\beta}$, then we can state convergence results for each
element in the support
$ \mathcal{B} = \{ u \in \Re: \Pr(\bar{\beta} = u) > 0 \}$ of
$\bar\beta \sim \pi_{\bar\beta}$. If $\beta \in \mathcal{B}$, h)
$s^{-1} \sum_{j=1}^p (\betady_j - \bbeta_{0,j}) \mathds{1}
\{\bbeta_{0,j} = \beta \} \overset{\textrm{a.s.}}{\longrightarrow}
(\mus - 1) \beta$, for
$\psi(t,u) = (t - (1 -\mus) u )\mathds{1} \{u = \beta \} /
\Pr\left(\bar{\beta} =~\beta\right)$; i)
$s^{-1} \sum_{j=1}^p (\betady_j / \mus - \bbeta_{0,j}) \mathds{1}
\{\bbeta_{0,j} = \beta \} \overset{\textrm{a.s.}}{\longrightarrow} 0$,
for
$\psi(t,u) = t\mathds{1} \{u = \beta \}/\{\mus \Pr\left(\bar{\beta}
  =~\beta\right)\}$, where $s = |\{j: \bbeta_{0,j} = \beta\}|$.
While the choices of $\psi(t,u)$ to derive h), i) are not
pseudo-Lipschitz, we can approximate $\mathds{1}\{\cdot\}$ by a smooth
function to arrive at the stated conclusion. If $\pi_{\bar\beta}$ is a
discrete-continuous distribution, then h) and i) hold for all $\beta$
in the support of the discrete
component. Figure~\ref{fig:supp-sc_scaled} of the Supplementary
Material document provides an additional illustration of the signal
recovery of $\betady / \mus$ similar to that of
Figure~\ref{fig:ra_scaled}. 

Note that, the existence of a solution $(\mus,\bs,\sigmas)$ to
\eqref{eq:state_evol} is required because the underlying AMP recursion
is specifically designed such that a fixed point corresponds to a
maximum of $\ell(\bbeta; \by^{*}, \bX)$; see
Section~\ref{subsub:gen_AMP} of the Supplementary Material
document.
That assumption can be verified numerically (see, for example, the
supplementary material of \citealp[Section~H.4.4]{sur+candes:2019} or
\citealp[Remark~2]{salehi+et+al:2019}), but proofs of the existence of stationary
points to state evolution equations have started appearing in recent works, such as
\citet{montanari+etal:2023}, \citet{li+sur:2024}, and
\citet{bellec+loriyama:2024}, and in settings that are similar in
spirit to ours. We expect the techniques used therein to be useful in
establishing the existence of a solution to
\eqref{eq:state_evol}. However, as we are able to obtain a numerical
solution $(\mus,\bs,\sigmas)$ for the wide range of
$(\alpha, \kappa,\gamma)$ settings that we have considered in our
extensive numerical studies, we defer these developments to future
work. The additional assumption that $\bb{J}(\mus,\bs,\sigmas)$ is
nonsingular arises from the smoothing argument that is necessary to
formally transfer asymptotic results from the AMP literature to
MDYPL. In particular, the limiting step at which we let our smooth
approximation converge to the AMP recursion of interest requires some
well-behavedness of the corresponding system of equations so that its
solution converges to $(\mus,\bs,\sigmas)$. We guarantee that this is
the case by requiring the nonsingularity of $\bb{J}(\mus,\bs,\sigmas)$
and using the Implicit Function Theorem. The assumption of
nonsingularity of $\bb{J}(\mus,\bs,\sigmas)$ is not present in
\citet{sur+candes:2019}, because the smoothing approximation step
has been missing in their analysis.

\subsection{Arbitrary covariate covariance}
\label{sec:arbitrary_covars}
We now relax the assumption of uncorrelated normal covariates and
instead assume that $\bX$ has rows that are realisations of
independent $\mathrm{N}(\b0_p,\bSigma)$ random vectors with
$\bSigma \in \Re^{p \times p}$. In what follows, let
$ \bSigma = \bb{L} \bb{L}^\top$, where $\bb{L}$ denotes the Cholesky
factor of $\bSigma$, and denote the minimum and maximum eigenvalues of
$\bSigma$ by $\lambda_{\min}(\bSigma)$ and $\lambda_{\max}(\bSigma)$,
respectively. As before, we require that
$ \var{(\bx_i^\top \bbeta_0)} = \bbeta_0^\top \bSigma \bbeta_0 \to
\gamma^2 \in \Re$ and that $p / n \to \kappa \in (0,1)$.

The MDYPL estimator inherits all equivariance properties the ML
estimator has. Hence, the derivations in \citet{zhao+etal:2022} enable
us to relate $\betady$ from covariates with covariance
$\bSigma = \bb{L}\bb{L}^\top$ and signal $\bbeta_0$ to the MDYPL
estimator $\thetady$ with $\mathrm{N}(\b0_p,\bb{I}_p)$ covariates and
signal $\btheta_0 = \bb{L}^\top \bbeta_0$ via
$\thetady = \bb{L}^\top \betady$.
In conjunction with Theorem~\ref{thm:aggregate_asymptotics}, we can
then get the following asymptotic result about the aggregate behaviour
of $\betady$.

\begin{theorem}
  \label{thm:aggregate_arbitrary}
  Consider the logistic regression model (\ref{eq:logistic_y}) with
  independent covariates $\bx_i \sim \mathrm{N}(\b0_p,\bSigma)$  and
  $\underset{n \to \infty}{\limsup } \, \lambda_{\max}(\bSigma) /
  \lambda_{\min}(\bSigma) < \infty$. Assume that
  $(\alpha,\kappa,\gamma)$ are such that \eqref{eq:state_evol} admits
  a solution $(\mus,\bs,\sigmas)$ such that
  $\bb{J}(\mus,\bs ,\sigmas)$ is nonsingular. Then, for any
  $t \in \Re$,
  \begin{equation}
    \frac{1}{p} \sum_{j=1}^p \mathds{1} \left\{	\sqrt{n} \tau_j \frac{\betady_{j}-\mus \bbeta_{0,j}}{\sigmas} \leq t \right\} \overset{p}{\longrightarrow} \Phi(t) \,,
  \end{equation}
  where $\Phi(\cdot)$ is the $\mathrm{N}(0,1)$ cumulative distribution
  function and $\tau_j^2 = \var\left(\bx_{ij} | \bx_{i, -j}\right)$
  is the conditional variance of the $j$th covariate given all others.
  In addition, if
  $\sum_{j=1}^p \delta_{\sqrt{n} \tau_j \bbeta_{0,j}} / p
  \overset{d}{\longrightarrow} \pi_{\bar{\beta}}$ and
  $n \sum_{j=1}^p \tau_j^2 \bbeta_{0,j}^2/ p
  \overset{p}{\longrightarrow} \expect(\bar{\beta}^2 )$
  for some distribution $\pi_{\bar{\beta}}$ with finite second moment
  and $\bar{\beta} \sim \pi_{\bar{\beta}}$, where $\delta_{(\cdot)}$
  is the Dirac delta function, then for any pseudo-Lipschitz function
  $\psi: \Re^2 \to \Re$ of order 2, 
  \begin{equation}
    \frac{1}{p} \sum_{j=1}^p  \psi \left(	\sqrt{n} \tau_j \left( \betady_{j}-\mus \bbeta_{0,j}\right), \sqrt{n} \tau_j \bbeta_{0,j} \right)  \overset{p}{\longrightarrow}  \expect \left[\psi(\sigmas G, \bar{\beta})\right], \quad \textrm{as } n \to \infty \,, 
  \end{equation}
  for $G \sim \mathrm{N}(0,1)$ independent of $\bar{\beta}$.
\end{theorem}

\subsection{Inference}
\label{subsec:inference}

The aggregate behaviour of $\betady$ in
Theorem~\ref{thm:aggregate_asymptotics} is the basis to extend the arguments
in \citet[Theorem~3.1 and Theorem~4.1]{zhao+etal:2022} and
\citet[Theorem~4]{sur+candes:2019} to MDYPL, and derive the asymptotic
distribution of adjusted $Z$-statistics and of DY prior penalized
likelihood ratio test statistics.

\begin{theorem}
  \label{thm:z_stats}
  Consider the logistic regression model~(\ref{eq:logistic_y}) with
  independent covariates $\bx_i \sim \mathrm{N}(\b0_p,\bSigma)$.
  Assume that $(\alpha,\kappa,\gamma)$ are such that
  \eqref{eq:state_evol} admits a solution $(\mus,\bs,\sigmas)$ such
  that $\bb{J}(\mus,\bs ,\sigmas)$ is nonsingular.  For any regression
  coefficient with $ \sqrt{n} \tau_j \bbeta_{0,j} = \mathcal{O}(1)$,
  $\sqrt{n} \tau_j (\betady_{j}-\mus \bbeta_{0,j}) / \sigmas
  \overset{d}{\longrightarrow} \mathrm{N}(0,1)$,
  where $\tau_j^2 = \var(\bx_{ij}|\bx_{i, -j})$.
\end{theorem}

Following the proof of \citet[Theorem~3.2]{zhao+etal:2022}, the
rotational invariance of normal random variables immediately extends
Theorem~\ref{thm:z_stats} to linear combinations of
parameters. Specifically, for any sequence of unit vectors
$\bv \in \Re^p$ such that
$\tau(\bv) \bv^\top \bbeta_0 = \mathcal{O}(n^{-1/2})$, it holds that $\sqrt{n} \tau(\bv) \bv^\top (\betady - \mus \bbeta_{0}) /\sigmas  \overset{d}{\longrightarrow} \mathrm{N}(0,1)$,
where
$\tau^2(\bv) = \var( \bv^\top \bx_i | \bb{P}_{\bv^\perp} \bv) =
(\bv^\top \bSigma^{-1} \bv)^{-1}$ and $\bb{P}_{\bv^\perp}$ is the
projection onto the subspace orthogonal to $\bv$. Most notably, this
extends Theorem~\ref{thm:z_stats} to arbitrarily large, but fixed,
collections of parameters. We now turn to the characterisation of the
DY prior penalized likelihood ratio (PLR) statistic.

\begin{theorem}
  \label{thm:likl_ratio_2}
  Consider the logistic regression model~(\ref{eq:logistic_y}) with
  independent covariates $\bx_i \sim \mathrm{N}(\b0_p,\bSigma)$. 
  Assume that $(\alpha,\kappa,\gamma)$ are such that
  \eqref{eq:state_evol} admits a solution $(\mus,\bs,\sigmas)$ such
  that $\bb{J}(\mus,\bs ,\sigmas)$ is nonsingular. Consider a fixed
  subset of indices $I = \{i_1,\ldots,i_k\}$, and the DY prior
  penalized likelihood ratio test statistic
  \begin{equation}
    \Lambda_I = \underset{\bbeta \in \Re^p}{\max} \, \ell(\bbeta;\by^*,\bX) - \underset{\substack{\bbeta \in \Re^p: \\ \bbeta_j = 0, \, j \in I }}{\max} \, \ell(\bbeta;\by^*,\bX) \,.
  \end{equation}
  Then, under the null that $\bbeta_{0,i_1} = \ldots = \bbeta_{0,i_k} = 0$, $\Lambda_I$ is asymptotically distributed as
  $2 \Lambda_{I} \overset{d}{\longrightarrow} \kappa \sigma_{*}^2/ b_{*} \chi^2_k$,
  where $\chi^2_k$ is a Chi-squared random variable with $k$ degrees of freedom. 
\end{theorem}

\begin{figure}
  \centering 
  \includegraphics[width = \linewidth]{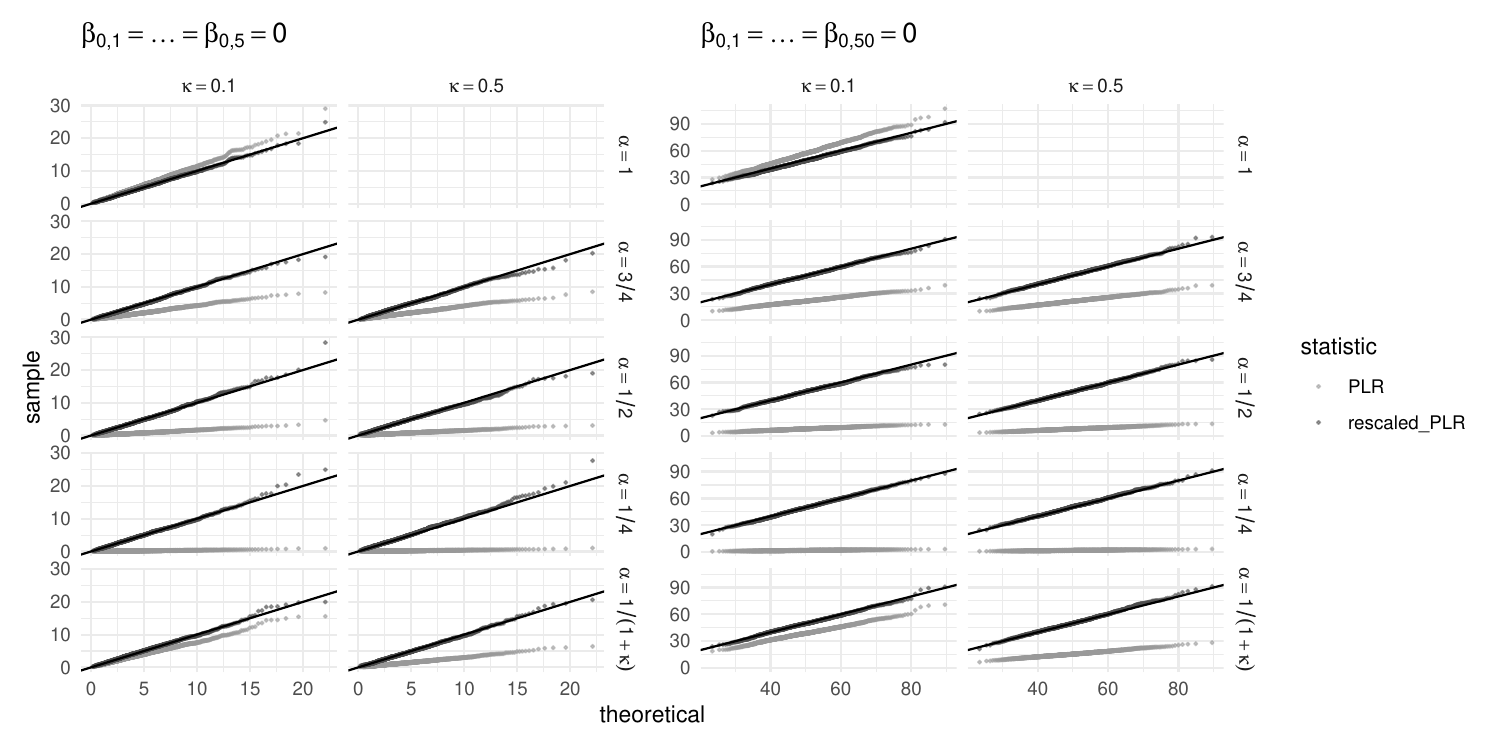}
  \caption{Q-Q plots comparing $\chi^2_k$ quantiles to empirical
    quantiles of the DY prior penalized likelihood ratio statistic
    (light grey) and its rescaled version
    $\bs / (\kappa \sigmas^2) \Lambda_{I}$ (dark grey), for
    $\alpha \in \{1, 3/4, 1/2, 1/4, 1/(1 + \kappa)\}$ and
    $I = \{1, 2, \ldots, 5\}$ ($k = 5$; left)
    $I = \{1, 2, \ldots, 50\}$ ($k = 50$; right).  The figures are
    based on $1000$ simulations of $\{\by, \bX\}$ where
    $\bx_i \sim \mathrm{N}(\b0_p,n^{-1}\bb{I}_p)$, $n = 2000$,
    $\kappa \in \{0.1, 0.5\}$, $p = n \kappa$, and $\bbeta_{0}$ has
    $p/2$ entries of zero and the remainder set to one, appropriately
    rescaled so that $\gamma^2 = 5$.}
  \label{fig:qqplots}
\end{figure}

We conduct a simulation study to assess the accuracy of the asymptotic
approximation for the distribution of the DY prior PLR statistic. We
generate $1000$ datasets $\{\by,\bX\}$, with $n = 2000$,
$\kappa \in \{0.1, 0.5\}$, $p = n\kappa$, covariate vectors
$\bx_i \sim \mathrm{N}(0,n^{-1} \bb{I}_{p})$, and signal $\bbeta_0$
having half of its entries equal to zero and the other half set to
one, appropriately rescaled so that $\gamma^2 = 5$. For each dataset,
we compute the DY prior PLR statistic using shrinkage parameter
$\alpha \in \{1, 3/4, 1/2, 1/4, 1 / (1 + \kappa)\}$ for the null
hypotheses $\bbeta_{0,1} = \ldots = \bbeta_{0,5} = 0$ and
$\bbeta_{0,1} = \bbeta_{0,2} = \ldots = \bbeta_{0,50} = 0$. According
to the results in \citet{candes+sur:2020}, the probability that the ML
estimator exists approaches one and zero for $\kappa = 0.1$ and
$\kappa = 0.5$, respectively. Hence, the methods in
\citet{candes+sur:2020}, i.e. when $\alpha = 1$, apply only for
$\kappa = 0.1$. Figure~\ref{fig:qqplots} shows the corresponding Q-Q
plots. As expected by Theorem~\ref{thm:likl_ratio_2}, we observe a
close agreement between the theoretical and empirical distributions
across all combinations of null hypotheses and shrinkage parameter
values.
	
\section{Shrinkage towards zero}
\label{sec:hyper}

\subsection{Adaptive shrinkage}
\label{subsec:adaptive_shrinkage}

\begin{figure}[t]
  \centering
  \includegraphics[width = \linewidth]{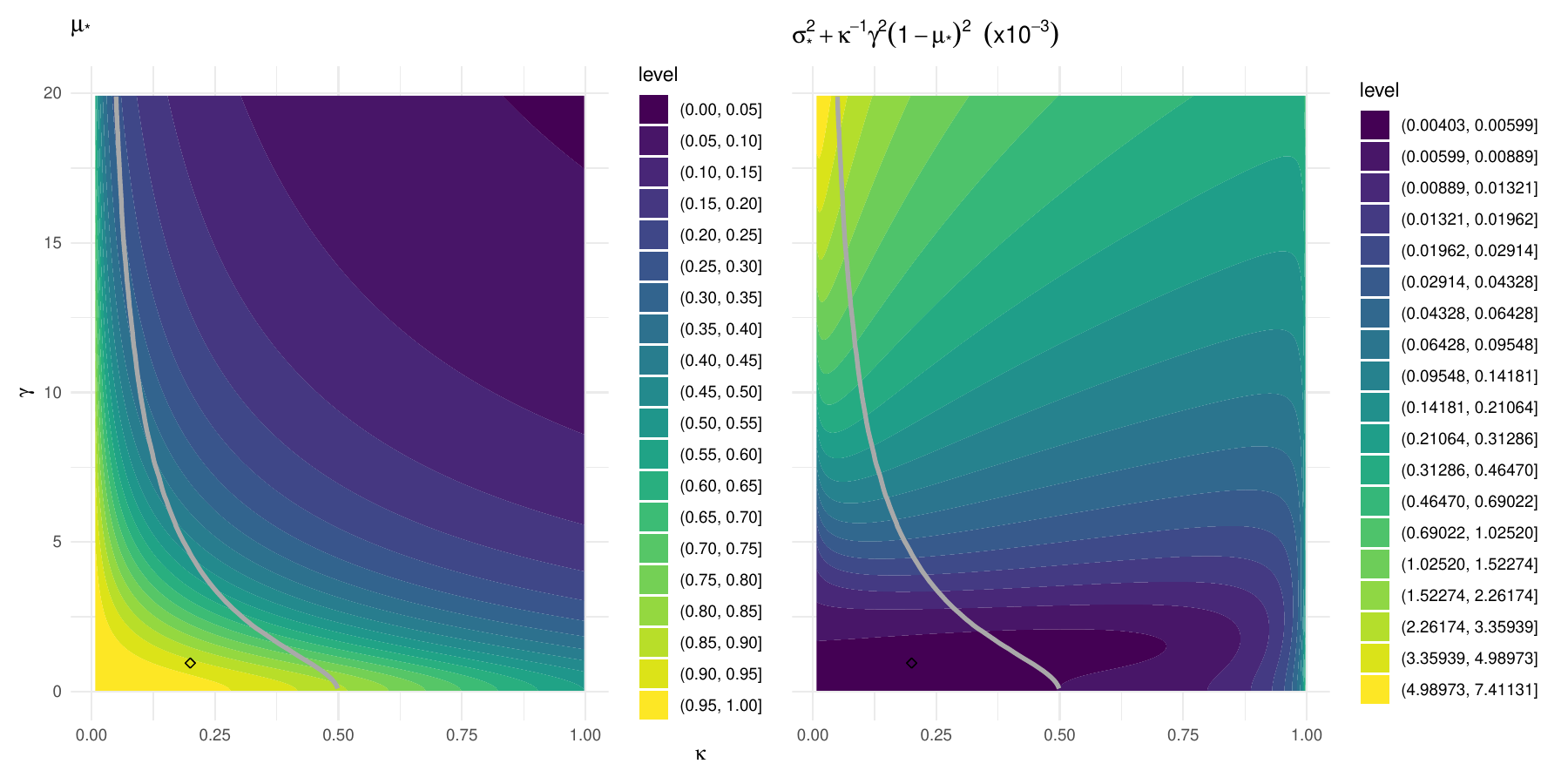}
  \caption{Contours of the asymptotic aggregate bias parameter $\mus$
    (left) and asymptotic MSE (right) for $\alpha = 1 / (1 + \kappa)$
    over a $(\kappa, \gamma)$ grid. The diamonds mark the point
    $(\kappa, \gamma) = (0.2, \sqrt{0.9})$, where the unscaled MDYPL estimator
    has been found to perform well in terms of aggregate bias
    ($\mu_{*} \approx 0.914$) and aggregate MSE
    ($\sigmas^2 + (1-\mus)^2 {\gamma^2}{\kappa^{-1}} \approx 5.08$) in
    the experiments of Section~\ref{sec:phase_transition}. The grey
    curve is the phase transition curve of \citet{candes+sur:2020}.}
  \label{fig:ra_pars}
\end{figure}

The MDYPL estimator not only exists for all data configurations
$\{\by,\bX\}$ and $\alpha \in (0, 1)$, but its dependence on $\alpha$
also allows one to tune the estimator to achieve desirable, or, in some
sense, optimal properties, by controlling the amount of shrinkage
induced as we vary $\alpha \in (0, 1)$.

Consider the shrinkage parameter $\alpha = 1 / ( 1+ \kappa)$, as
suggested in \citet{rigon+aliverti:2023}. This choice of $\alpha$ is
adaptive in that the larger the limiting proportion $\kappa$ of
covariates to observations, the more the MDYPL estimates shrink
towards zero.

We revisit the setting of empirical study in
Section~\ref{sec:phase_transition}. Figure~\ref{fig:ra_pars} shows
surface plots of the asymptotic aggregate bias parameter $\mus$ as
well as the asymptotic aggregate MSE (aMSE)
$\sigmas^2 + (1-\mus)^2 {\gamma^2}{\kappa^{-1}}$ in the
$(\kappa,\gamma)$ plane for $\alpha = 1 / (1+ \kappa)$. At
$\kappa = 0.2$ and $\gamma = \sqrt{0.9}$, the asymptotic aggregate
bias parameter and root aMSE are $\mus \approx 0.91$ and
$\sqrt{\sigmas^2 + (1-\mus)^2 {\gamma^2}{\kappa^{-1}}} \approx 2.25$,
respectively, which are close to $\mu^* = 1$ (aggregate unbiasedness) and the
lower bound $\sigma_{*} \approx 2.22$ for the root aMSE of the MDYPL
estimator.

As observed in Section~\ref{sec:phase_transition}, the favourable
properties of the MDYPL estimator, namely low bias and low aMSE at
$(\kappa, \gamma) = (0.2, \sqrt{0.9})$ do not generalise to other
$(\kappa, \gamma)$ values. Indeed, Figure~\ref{fig:ra_scaled} shows
that $\mus$ takes smaller and smaller values as $\kappa$ and $\gamma$
increase, and, hence, the shrinkage induced by
$\alpha = 1 /(1 + \kappa)$ is
excessive. Section~\ref{subsec:default_mdypl_summaries} of the
Supplementary Material document provides further evidence in that
direction.

Overall, setting $\alpha = 1/(1 + \kappa)$ is attractive because we
are able to find solutions to equations \eqref{eq:state_evol} for a
wide range of $(\kappa, \gamma)$ values. However, it is inadequate as
a means of bias reduction in the high-dimensional setting without the
rescaling of the MDYPL estimator predicted by
Theorem~\ref{thm:aggregate_asymptotics}.

\subsection{Aggregate unbiasedness}
\label{subsec:bias}

\begin{figure}[t]
  \centering
  \includegraphics[width = \textwidth]{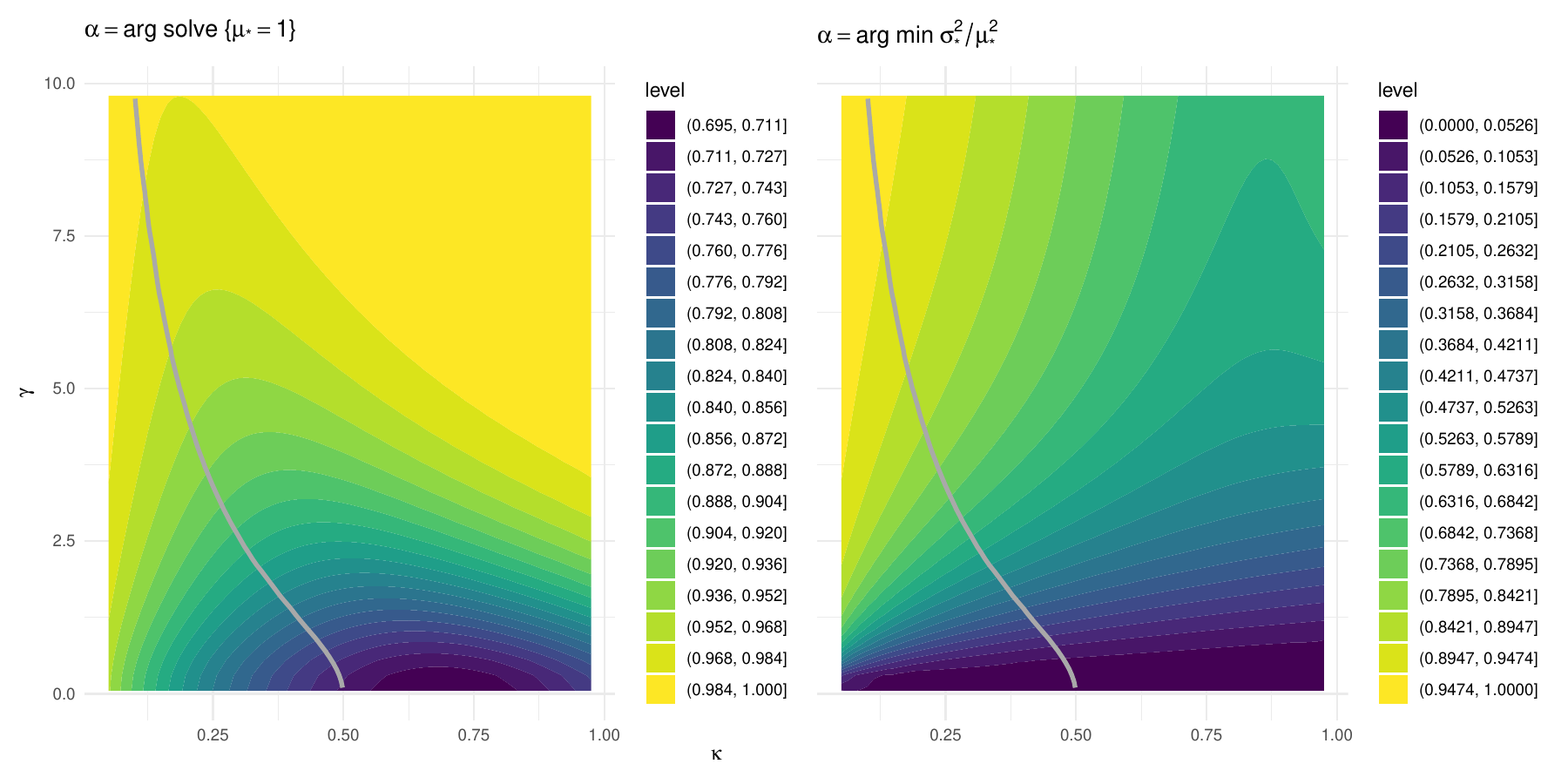}
  \caption{Contours of the values of shrinkage parameter $\alpha$ that
    achieve asymptotic aggregate unbiasedness $\mus = 1$ (left), and
    minimal asymptotic aggregate variance ${\sigmas^2} / {\mus^2}$ of
    $\betady / \mus - \bbeta_0$ (right). The grey curve is the phase
    transition curve of \citet{candes+sur:2020}.}
  \label{fig:opt_pars}
\end{figure}  

As we observed, for specific regions of $(\kappa, \gamma)$, the
scaling $\alpha = 1/(1 + \kappa)$ produces estimates that are close to
being asymptotically unbiased in an aggregate sense
$(\mus \approx 1)$. This naturally raises the question what scaling
$\alpha$ produces estimates that have zero aggregate bias. We
investigate that by numerically searching for $\alpha$ values such
that $(1, \bs, \sigmas)$ is a solution to \eqref{eq:state_evol}.

The left panel of Figure~\ref{fig:opt_pars} shows a contour plot of
$\alpha$ for which $\mus = 1$ on the $(\kappa, \gamma)$ plane. We note
that as $\kappa$ increases, a higher amount of shrinkage (smaller
$\alpha$) is required to achieve unbiased estimates. On the other
hand, as $\gamma$ grows, less shrinkage ($\alpha \to 1$) is required
to achieve unbiased estimates. These observations are intuitive; the
bias of the MDYPL estimator increases with $\kappa$, and decreases
with the signal strength.

\subsection{Minimal mean squared error}
\label{subsec:var} 

Theorem~\ref{thm:aggregate_asymptotics} suggests that we can eliminate
the aggregate asymptotic bias of the MDYPL estimator for any given
$\alpha \in (0, 1)$ by rescaling it. Hence, another natural question
is what value of $\alpha$ minimises the limit $\sigmas^2/\mus^2$ of
the aMSE $\vnorm{\betady / \mus - \bbeta_0}_2^2/ p$ of the rescaled
estimator.

Towards answering this question, we numerically solve \eqref{eq:state_evol} for
all combinations of $\kappa \in \{0.05, 0.0625, \ldots, 0.975\}$,
$\gamma \in \{0.05, 0.3, \ldots, 10\}$, and $100$ values for $\alpha$
ranging from $0.01$ to $0.99$.  For each $(\kappa, \gamma)$
combination, we identify the $\alpha$ value that minimises
$\sigmas / \mus$ among the $\alpha$ values. The resulting contour plot
is shown in the right panel of Figure~\ref{fig:opt_pars}. As $\kappa$
tends to zero (and consequently $\mus \to 1$), the $\alpha$ that
minimises the asymptotic MSE of the rescaled MDYPL estimator
approaches one, hinting at the fact that in a low-dimensional setting,
the ML estimator achieves asymptotically optimal variance.

\section{Comparison to other methods}
\label{sec:other_methods}

\subsection{Preamble}

In this section, we compare the performance of the MDYPL estimator
with the corrected least-squares (CLS) estimator of
\citet{lewis+battey:2023} and the logistic ridge of
\citet{salehi+et+al:2019} whose asymptotic properties have also been
studied in the $p/n \to \kappa \in (0,1)$ logistic regression
setting. For the comparison we use simulation settings that have been
used in those works.

\subsection{Corrected least-squares estimator}
\label{subsec:corrected_ls}

The CLS estimator of \citet{lewis+battey:2023} relies on a consistent
estimator of $\bb{\eta}_0 = \bX
\bbeta_{0}$. Leaving technical details to \citet{lewis+battey:2023},
the CLS estimator is of the form
$\hat{\bbeta}^{\ttextrm{CLS}} = \varsigma^{-1}
\left(\hat{\bbeta}^{\ttextrm{LS}} -\bb{\delta} \right)$, where
$\hat{\bbeta}^{\ttextrm{LS}} = (\bX^\top \bX)^{-1} \bX^\top (2\by -
{\bb 1}_n)$, where ${\bb 1}_n$ is a vector of $n$ ones, and
$\varsigma, \bb{\delta}$ are scale and location corrections to
$\hat{\bbeta}^{\ttextrm{LS}}$.  In the fixed design case, and under
some regularity conditions,
\citet[][Sections~5.2-5.3]{lewis+battey:2023} establish consistency of
the CLS estimator in $\ell_\infty$ and $p^{-1/2}$-scaled $\ell_2$
norm.  \citet[][Section~5.4 and Section~8]{lewis+battey:2023} further
conjecture that the test statistic
\begin{equation}
  \label{eq:cls_t_stat}
  T_j = \frac{	\hat{\bbeta}^{\ttextrm{CLS}}_j - \bbeta_{0, j}  }{\varsigma^{-1} \vnorm{\bb{e}_j^\top (\bX^\top \bX)^{-1} \bX \bb{\Gamma}^{1/2}}_2}
\end{equation}
is asymptotically distributed according to $\mathrm{N}(0,1)$, where
$\bb{e}_j$ is a vector of zeros apart from its $j$th element that is
one, and $\bb{\Gamma}$ is the variance-covariance matrix of
$2\by - {\bb 1}_n$.

Fundamental to those asymptotic statements are assumptions about the
data generating process and the use of a consistent estimator for the
linear predictor $\bb{\eta}_0 = \bX \bbeta_{0}$ according to a certain
divergence \citep[see][Conditions 1-2]{lewis+battey:2023}, in order to
define appropriate estimators for $\varsigma$, $\bb{\delta}$. For the
random design case of Section~\ref{sec:arbitrary_covars},
\citet[Section~8]{lewis+battey:2023} provide results for estimating
$\bb{\eta}_0$ using the logistic LASSO, which maximises
$\ell(\bbeta; \by, \bX) - \lambda \sum_{j=1}^p | \bbeta_j | $. In
particular, if C1) $\lambda_{\min}(\bSigma) > c_1$ and
$\max_{1 \leq j \leq p} \bSigma_{jj} < c_2 $ for some
constants $c_1, c_2 > 0$, and C2) $\vnorm{\bbeta_{0}}_\infty = \mathcal{O}(1)$,
  $\vnorm{\bbeta_{0}}_0 = \mathcal{O}(n^{1/2 - c_3})$ for some
  constant $c_3 > 0$ and where $\vnorm{\bbeta_{0}}_0 $ is the support
  of $\bbeta_{0}$,
then the LASSO CLS achieves consistency.

We note that condition C1) is similar in spirit to
${\lim \sup}_{n \to \infty} \, \lambda_{\max}(\bSigma) /
\lambda_{\min}(\bSigma) < \infty$, which is required for
Theorem~\ref{thm:aggregate_arbitrary}. The first part of condition
C2) is in line with our requirement that the empirical
distribution of $\bbeta_{0}$ converges to some distribution
$\pi_{\bar{\beta}}$. The sparsity assumption
$\vnorm{\bbeta_{0}}_0 = \mathcal{O}(n^{1/2 - c_3})$ of C2),
which is required to achieve consistent estimation of $\bb{\eta}_0$
using LASSO, is stronger than what is required for MDYPL, which is
agnostic about the sparsity of $\bbeta_{0}$ so that $\bbeta_{0}$ can even
have full support. From a practical standpoint, the choice of the
regularisation parameter $\lambda$ for the LASSO, that achieves
consistent estimates of the linear predictors, may be difficult. For
instance, \citet[Proposition~7]{lewis+battey:2023} give conditions for
consistent estimation of $\bb \eta_0$ with
$\lambda = A \left\{(\log p \log n) / n \right\}^{1/2}$, for some
constant $A$ that remains unspecified.

To compare the performance of MDYPL to CLS, we consider a setting
close to the experiment in \citet[Supplementary Material,
Section~4.1]{lewis+battey:2023}. We set $\gamma = 3$, and consider
$\kappa \in \{0.2, 0.5 \}$ and $n \in \{400, 800, \ldots, 2000 \}$,
and set $p = \kappa n$ for each combination of $\kappa$ and $n$.  For
each setting of $(n, p)$, we draw the model matrix $\bX$ with rows
$\bx_i \sim \mathrm{N}(\b0_p, \bSigma)$, where $\bSigma$ has diagonal
elements $\bSigma_{jj} = 1$ and off-diagonal elements
$\bSigma_{jk} = 0.5$ $(j \neq k)$. The true parameter vector
$\bbeta_0$ has all but its first five elements set to zero. The
nonzero elements are all equal and scaled such that
$\bbeta_{0}^\top \bSigma \bbeta_{0} = \gamma^2 = 9$. For each
combination of $\kappa$ and $n$, we draw $5000$ independent samples of
responses $\by$ according to \eqref{eq:logistic_y}. We evaluate the
performance of the CLS with oracle linear predictors
$\bb{\eta}_0 = \bX \bbeta_0$ and with estimated linear prediction
coming from LASSO. For the LASSO, we used the R \citep{R} package
\texttt{glmnet} \citet{friedman+etal:2015} with default options for
estimation and selection of the turning parameter.  For the rescaled
MDYPL estimator, we set $\alpha = 0.95$, and consider its performance
with oracle parameters $\mus, \bs, \sigmas, \tau_j$, which require
knowledge of $\gamma$ and $\bSigma$. Note that the
oracle CLS estimator has knowledge of the true linear predictor
$\bb{\eta}_0$, while the oracle MDYPL only requires knowledge of
$\kappa$, $\gamma$, and $\tau_j$. In this setting,
$\tau_j^2 = (p + 1) / 2p$ for all $j$ and, thus, by
Theorem~\ref{thm:aggregate_arbitrary},
\begin{equation}
  \frac{n + \kappa^{-1}}{2} \frac{\vnorm{\betady / \mus - \bbeta_{0}}_2^2}{p} \overset{p}{\longrightarrow} \frac{\sigmas^2}{\mus^2}, \quad \textrm{as } n \to \infty \,.
\end{equation} 
As a result, the rescaled MDYPL estimator also achieves the same $p^{-1/2}$-scaled
$\ell_2$-norm consistency as the LASSO CLS and oracle CLS estimators
do, without exploiting sparsity as an assumption. We also consider the
non-oracle rescaled MDYPL estimator where we estimate the unknown
constants $\tau_j$ and
$\eta = \sqrt{\mus^2 \gamma^2 + \kappa \sigmas^2}$ from the data using
the procedures we introduce in Section~\ref{sec:unknown_constants}.

\begin{figure}[t]
  \centering
  \includegraphics[width=\textwidth]{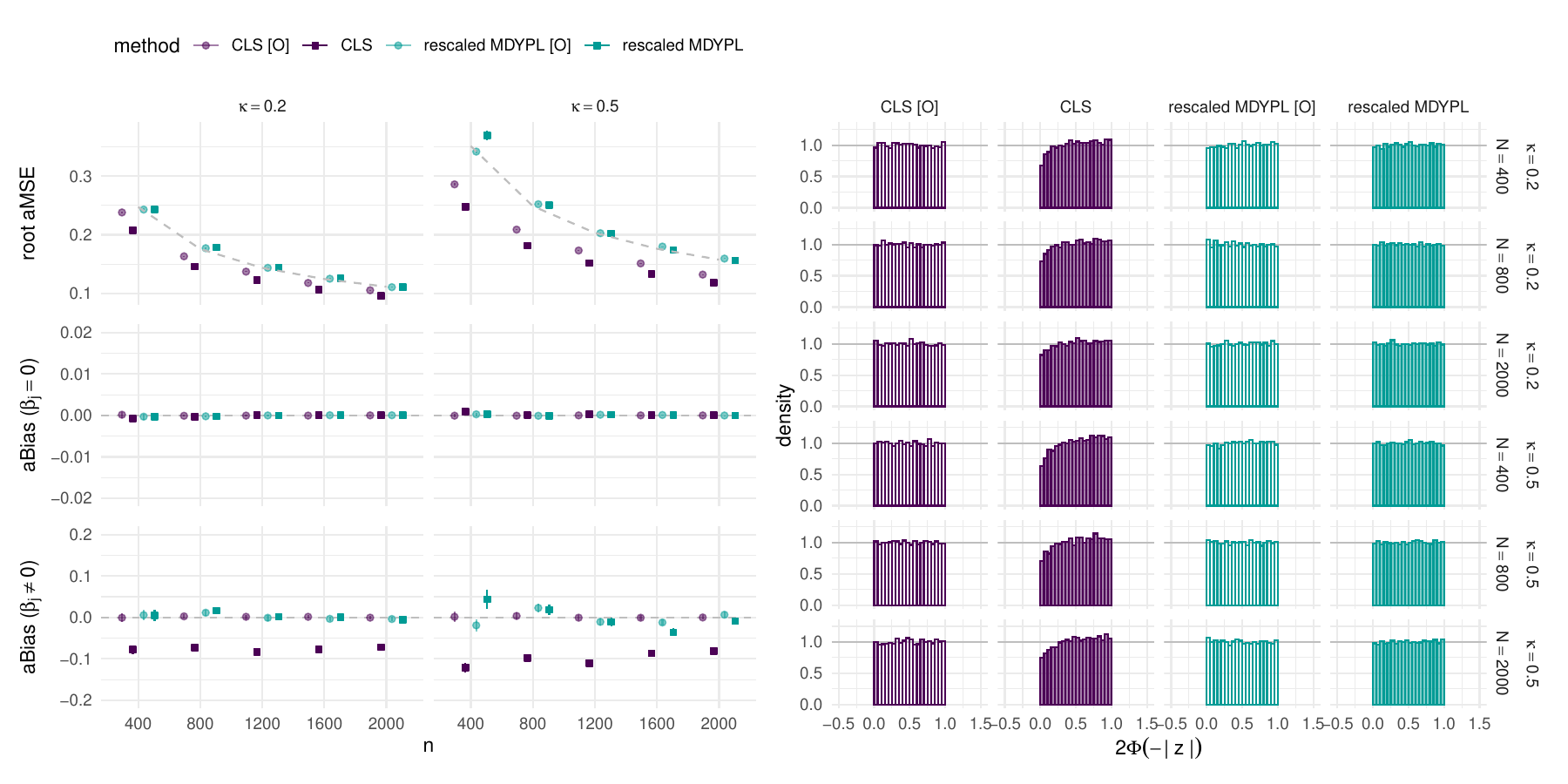}
  \caption{Performance comparison of the oracle and non-oracle
    versions of the CLS and MDYPL procedures for estimation and
    inference in the simulation setting of
    Section~\ref{subsec:corrected_ls}.  The top left panel shows the
    estimated root aggregate MSE (root aMSE) of the estimators
    for $n \in \{400, 800, \ldots, 2000\}$, $\kappa \in \{0.2, 0.5\}$. The grey dashed curve represents the asymptotic root aMSE
    $\sigmas \mus^{-1} \sqrt{2 (n + \kappa^{-1})^{-1}}$ of the rescaled
    MDYPL estimator. The mid-left and bottom-left panels show the
    aggregate bias (aBias) for the zero and non-zero elements of the
    parameter vector $\bbeta_{0}$, respectively. The right panel shows
    the estimated finite-sample distributions of the $p$-values from
    the two-sided test that a single element of the parameter vector
    $\bbeta_{0}$ is zero, for the first $25$ zero parameters.}
  \label{fig:cls_signal}
\end{figure}

Figure~\ref{fig:cls_signal} summarises the performance of the
considered estimators for signal recovery and inference. On the left
panel of Figure~\ref{fig:cls_signal}, we see that CLS outperforms
rescaled MDYPL in root aMSE, and that the estimated finite-sample root
aMSE of MDYPL closely agrees with its asymptotic limit. In contrast to
CLS, the performance of the oracle and non-oracle MDYPL estimator is
effectively the same. Furthermore, all methods deliver almost unbiased
estimators of the zero elements of the parameter vector. However, we
see that the low estimation error of the LASSO CLS estimator is
accompanied by biased estimates of the non-zero elements of the
parameter vector. That bias appears to decrease slower than the usual
$O(n^{-1})$ asymptotic rate. On the right panel of
Figure~\ref{fig:cls_signal}, we see the estimated finite-sample
distributions of the $p$-values from the two-sided test that a single
element of the parameter vector $\bbeta_{0}$ is zero, for the first
$25$ zero parameters, based on the oracle and non-oracle versions of
the CLS-based statistic~(\ref{eq:cls_t_stat}) and of the MDYPL-based
statistic in Theorem~(\ref{thm:z_stats}). We observe that the oracle
version of the CLS-based statistic produces almost uniform $p$-value
distributions across the considered $(\kappa, n)$ settings. However,
deviations are observed for the non-oracle version in the direction of
producing conservative inferences, with the distributions appearing to
converge only slowly to uniform as $n$ increases. In stark contrast,
both oracle and non-oracle versions of the MDYPL-based statistic
result in almost uniform $p$-value distributions across the
considered $(\kappa, n)$ settings. 

\subsection{Logistic ridge}
\label{subsec:reg_lr}

\begin{figure}[t]
  \centering 
  \includegraphics[width = \linewidth]{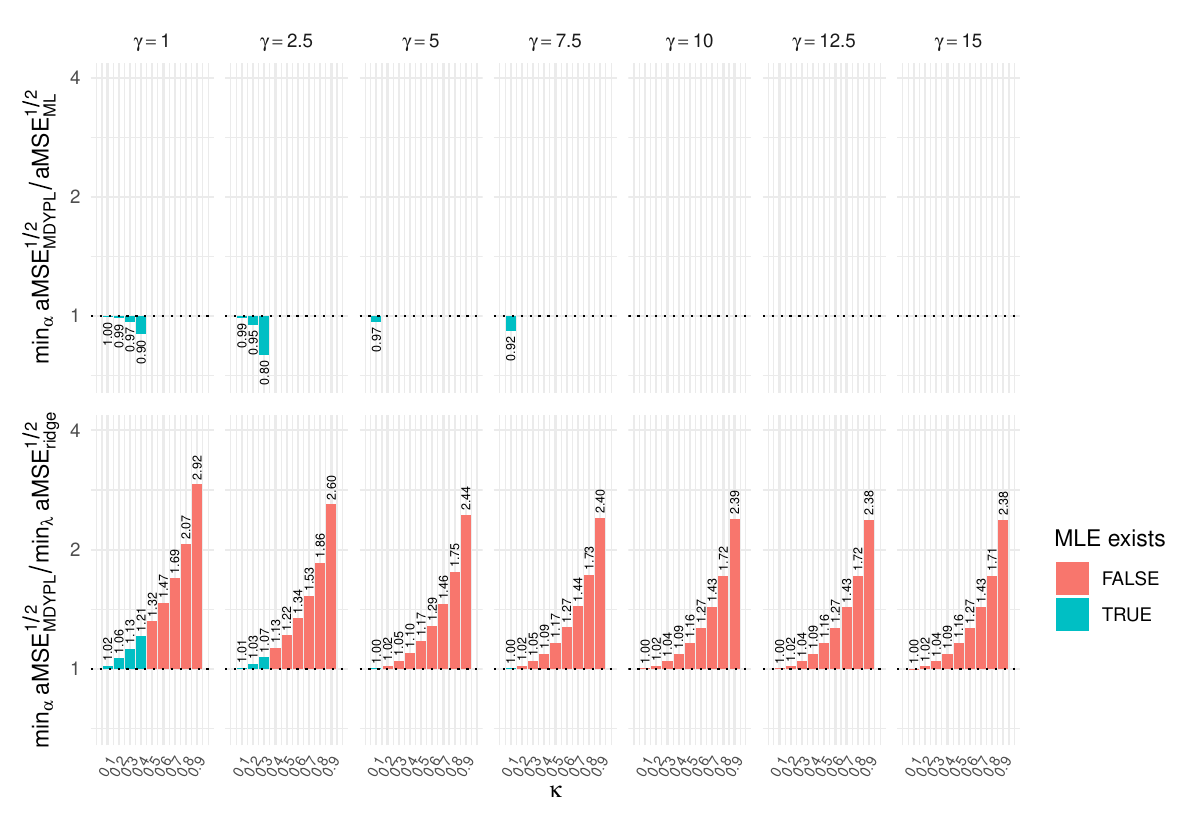}
  \caption{The square root of the minimum asymptotic aggregate MSE
    (aMSE) of the rescaled MDYPL estimator $\betady(\alpha) / \mus$
    over $\alpha$, relative to that of the rescaled ridge estimator
    $\betar(\lambda) / \bar{\mu}$ over $\lambda$ (bottom), and of the
    asymptotic aMSE of the rescaled ML estimator $\betam / \mus$
    (top). The $y$-axes are in log-scale. We consider all combinations
    of $\kappa \in \{0.1, 0.2, \ldots, 0.9\}$ and
    $\gamma \in \{1, 2.5, 5, 7.5, \ldots, 15\}$. For each
    $(\kappa, \gamma)$ setting, the asymptotic aMSE for MDYPL and
    logistic ridge regression are computed over grids of $500$ values
    for $\alpha \in (0, 1)$ and $\lambda \in (0, 4.5)$,
    respectively. See Section~\ref{subsec:reg_lr} for details.}
  \label{fig:rlr_mse_scaled}
\end{figure}

\citet[Chapter~4]{sur:2019} provides insights into the behaviour of
the ridge estimator in $p / n \to \kappa \in (0, 1)$ logistic
regression. Later, \citet{salehi+et+al:2019} formalised the aggregate
asymptotic behaviour of ridge regression for
$p/n \to \kappa \in (0,\infty)$, with a corresponding result as that
in Theorem~\ref{thm:aggregate_asymptotics}. As is the case for the
MDYPL estimator, the amount of shrinkage that is induced by ridge
regression is controlled by a single tuning parameter. So, it is of
interest to compare the performance of the MDYPL and ridge regression
estimators.

We use the setting in \citet{salehi+et+al:2019}. Let $\{\by,\bX\}$ be
data from a logistic regression model with covariates
$\bx_i \sim \mathrm{N}(\b0_p,p^{-1}\bb{I}_p)$ and signal $\bbeta_0$,
the entries of which are realizations of independent random variables
with distribution $\pi_{\bar{\beta}}$, such that the signal strength
$\vnorm{\bbeta_0}_2^2 / p$ is equal to $ \gamma^2$, as
$p / n \to \kappa \in (0,1)$. The ridge logistic regression estimator
$\betar$ is the maximiser of $\ell(\bbeta; \by, \bX) / n - \lambda \sum_{j=1}^p \bbeta_{j}^2 / (2p)$,
where $\lambda \geq 0$ is the tuning parameter that controls the
amount of shrinkage induced by the $\ell_2$
penalty. \citet[Theorem~2]{salehi+et+al:2019} show that for any
locally Lipschitz function $\psi:\Re^2 \to \Re$,
\begin{equation}
  \label{eq:ridge_asymptotics}
  \frac{1}{p} \sum_{j = 1}^p \psi \left(\betar_j - \bar{\mu} \bbeta_{0,j}, \bbeta_{0,j} \right) \overset{p}{\longrightarrow} \expect \left[\psi \left(\bar{\sigma} G, \bar{\beta} \right)\right] \,, 
\end{equation} 
where $G \sim \mathrm{N}(0,1)$, $\bar{\beta} \sim \pi_{\bar{\beta}}$
independent of $Z$ and $\bar{\mu},\bar{\sigma}$ are part of the
solution $(\bar{\mu}, \bar{b}, \bar{\sigma})$ to the system of
nonlinear equations in $(\mu,b,\sigma)$
\begin{equation} 
  \begin{aligned}
    \label{eq:ridge_FOC_}
    \expect\left[2 \zeta'' \left( - \gamma Z_1 \right) \prox{b \zeta}{{\mu}\gamma Z_1 + {\sigma}Z_2} \right] + {\mu} \kappa &= 0\\
    1 - \kappa + {b} \lambda - \expect \left[\frac{2 \zeta'\left(-\gamma Z_1\right)}{1+ {b} \zeta'' \left(\prox{{b} \zeta}{{\mu}\gamma Z_1 + {\sigma}Z_2}\right)} \right]  &= 0 \\ 
    {\sigma}^2 \kappa - \expect \left[2 \zeta' \left(- \gamma Z_1 \right) \left\{{\mu}\gamma Z_1 + {\sigma}Z_2 - \prox{{b} \zeta}{{\mu}\gamma Z_1 + {\sigma}Z_2}\right\}^2 \right]  &= 0 \,,  
  \end{aligned}
\end{equation}
where $Z_1$ and $Z_2$ are independent $\mathrm{N}(0,1)$ random
variables. As is the case for $\betady$ with $\alpha = 1$,
$\betar$ is simply the ML estimator for $\lambda = 0$. In fact, it can
be shown that equations~(\ref{eq:ridge_FOC_}) with $\lambda = 0$
coincide with \eqref{eq:state_evol} with $\alpha = 1$ and
$(\mus,\bs,\sigmas) = (\bar{\mu},\bar{b}, \kappa^{-1/2} \bar{\sigma}
)$. Lemma~\ref{lemma:equivar} of the Supplementary Material document
and Theorem~\ref{thm:aggregate_asymptotics} establish that, for any
pseudo-Lipschitz function $\psi: \Re^2 \to \Re$ of order two, the
MDYPL estimator behaves as
\begin{equation}
  \label{eq:dy_ridge}
  \frac{1}{p} \sum_{j=1}^{p} \psi\left(\kappa^{-1/2} \left(\betady_j - \mus \bbeta_{0,j}\right), \kappa^{-1/2}\bbeta_{0,j} \right) \overset{\textrm{a.s.}}{\longrightarrow} \expect \left[\psi(\sigmas G, \kappa^{-1/2} \bar{\beta} )\right], \quad \textrm{as } n \to \infty \,,
\end{equation}
where $(\mus,\sigmas)$
solve \eqref{eq:state_evol}, and $G \sim \mathrm{N}(0,1)$ is
independent of $\bar{\beta} \sim \pi_{\bar{\beta}}$.

Using \eqref{eq:ridge_asymptotics} and \eqref{eq:dy_ridge}, we can
compare $\betady$ and $\betar$ based on various performance
metrics. Here, we compare the minimum asymptotic aMSEs of the scaled
estimators $\betady/\mus$ and $\betar / \bar{\mu}$, which are
$\min_{\alpha \in (0, 1)} \kappa \sigmas^2 / \mus^2$ and
$\min_{\lambda \in (0, \infty)}\bar{\sigma}^2 / \bar{\mu}^2$,
respectively.

We consider all combinations of
$\gamma \in \{1, 2.5, 5, 7.5, \ldots, 15\}$ and
$\kappa \in \{0.1, 0.2, \ldots, 0.9\}$. For each $(\kappa, \gamma)$
specification, the asymptotic aMSE is computed over grids of $500$
values for $\alpha \in (0, 1)$ and $\lambda \in (0, 4.5)$.  We were
able to obtain solutions $(\mus, \bs, \sigmas)$ and
$(\bar{\mu}, \bar{b}, \bar{\sigma})$ across all settings we
considered. Also, none of the $\alpha$ and $\lambda$ values that
minimise the asymptotic aMSEs over the respective grids are at the
endpoints of the grids.

Figure~\ref{fig:rlr_mse_scaled} shows the square root of the minimum
asymptotic aMSE of the rescaled MDYPL estimator relative to that of
the rescaled ridge estimator, and relative to the root of the
asymptotic aMSE of the rescaled ML estimator, whenever the latter
exists. Values smaller than $1$ indicate that the rescaled MDYPL
estimator performs better in terms of minimum asymptotic aMSE.

The rescaled MDYPL estimator appears to perform similarly to the
rescaled ML estimator, whenever the latter exists, for small $\kappa$
and small $\gamma$ values, but substantially outperforms it as
$\kappa$ and $\gamma$ increase. On the other hand, the rescaled MDYPL
estimator performs similarly to the rescaled ridge estimator for small
to moderate $\kappa$, but the latter increasingly outperforms MDYPL
for large values of $\kappa$. The ratio in minimum asymptotic aMSEs
seems to converge to fixed values as $\gamma$ increases. Logistic
ridge regression results in better asymptotic aMSE performance than
MDYPL, most likely due to the stronger shrinkage it imposes towards
zero. To see that, ridge regression can be thought of as adding $p$
new observations to each step of the re-weighted least-squares
procedure for solving the ML equations, which is a considerable amount
of information in proportional asymptotics and can return estimates
even for $\kappa > 1$. On the other hand, MDYPL simply adjusts the
existing $n$ responses by a small constant. Overall, we find that the
rescaled MDYPL estimator is able to maintain a minimum MSE that is
comparable to that of the rescaled ridge estimator, unless $\kappa$ is
large. This is important given the availability of the asymptotically
valid procedures for inference in Theorem~\ref{thm:likl_ratio_2} for
MDYPL, which, to the knowledge of the authors, are not available for
logistic ridge regression. 

\section{Inclusion of an intercept}
\label{sec:intercept}

In most practical applications of logistic regression, the inclusion
of an intercept parameter improves model fit and interpretability. The
logistic regression model with an intercept parameter has
\begin{equation}
  \label{eq:logistic_intercept}
  \Pr(y_i = 1 \mid \bx_i) = \zeta'(\theta_0 + \bx_i^\top \bbeta_{0}) \quad (i = 1, \ldots, n) \,,
\end{equation}
for $\theta_0 \in \Re$, and is not immediately covered by the
developments in the previous sections. We present a conjecture about
the aggregate behaviour of the MDYPL estimator in the presence of an
intercept term in the linear predictor. The conjecture can be obtained
mechanistically by contrasting the state evolution equations for the
MDYPL estimator in expression~(\ref{eq:state_evol}) to those of
\citet[expression~(5)]{sur+candes:2019} for the ML estimator, and use
the found relationships to extend \citet[Conjecture
7.1]{zhao+etal:2022} about the behaviour of the ML estimator with
intercept to the case of the MDYPL estimator. The resulting system of
nonlinear equation in $(\mu,b,\sigma,\iota)$ is
\begin{equation}
  \label{eq:state_evol_inter}
  \begin{aligned}
    0 & = \expect \left[\zeta'(Z_1)Z_1 Q_{+}(\alpha, b, Z_2)  \right] - \expect \left[\zeta'(-Z_1)Z_1 Q_{-}(\alpha, b, Z_2)  \right] \\ 
    1 - \kappa &= \expect \left[ \frac{\zeta'(Z_1)}{1 + b \zeta'' \left( \prox{b \zeta}{\frac{1 + \alpha}{2}b + Z_2 } \right)}  \right] + \expect \left[ \frac{\zeta'(-Z_1)}{1 + b \zeta'' \left( \prox{b \zeta}{\frac{1 + \alpha}{2} b - Z_2 } \right)}  \right] \\ 
    \frac{\sigma^2 \kappa^2}{b^2} &= \expect \left[\zeta'(Z_1) Q_{+}(\alpha, b, Z_2)^2  \right] + \expect \left[\zeta'(-Z_1) Q_{-}(\alpha, b, Z_2)^2  \right] \\ 
    0 &= \expect \left[ \zeta'(Z_1) Q_{+}(\alpha, b, Z_2) \right] - \expect \left[ \zeta'(-Z_1) Q_{-}(\alpha, b, Z_2) \right] \,, 
  \end{aligned}
\end{equation}
with
\[
  Q_{+}(\alpha, b, Z) = \frac{1 + \alpha}{2} - \zeta' \left( \prox{b \zeta}{\frac{1 + \alpha}{2} b + Z} \right) \,, \quad
  Q_{-}(\alpha, b, Z) = \frac{1 + \alpha}{2} - \zeta' \left( \prox{b \zeta}{\frac{1 + \alpha}{2}b - Z} \right) \,,
\]
and $Z_1 = \gamma Z + \theta_0$,
$Z_2 = \mu \gamma Z + \sqrt{\kappa} \sigma G + \iota$, where $Z$ and $G$ are
independent $\mathrm{N}(0,1)$ random variables, and $\theta_0$
is the intercept term in the linear predictor.

\begin{conjecture}
  \label{conj:intercept}
  Consider the logistic regression model with intercept in
  \eqref{eq:logistic_intercept}, and where covariates are independent
  realisations $\bx_i \sim \mathrm{N}(\bb{0}_p, \bSigma)$ with
  $\bbeta_{0}, \bSigma$ such that
  $\bbeta_{0}^\top \bSigma \bbeta_{0} \to \gamma^2$, as
  $n \to \infty$.
  Let $(\hat{\theta}_0^{\ttextrm{DY}}, \betady)^\top$ be the MDYPL
  estimator of a logistic regression model with intercept, where
  $\hat{\theta}_0^{\ttextrm{DY}}$ and $\betady$ denote the MDYPL
  estimates of the intercept and all other coordinates,
  respectively. Assume that $(\alpha, \kappa, \gamma, \theta_0)$ are
  such that a nonsingular solution $(\mus, \bs, \sigmas, \iota_{*})$
  to \eqref{eq:state_evol_inter} exists. Then, $\hat{\theta}_0^{\ttextrm{DY}} \overset{p}{\longrightarrow} \iota_{*}$ as $n \to \infty$; and, Theorem~\ref{thm:z_stats} and Theorem~\ref{thm:likl_ratio_2} hold with $\mus, \bs, \sigmas$ coming from \eqref{eq:state_evol_inter}.
\end{conjecture}

\begin{table}[t]
  \caption{Estimates of
    $P(b_{*} 2 \Lambda_{I} / (\kappa \sigma_*^2) \le \chi^2_{10, q})$
    (in percentage) from the simulation setting of
    Section~\ref{sec:intercept}, where $\Lambda_{I}$ is the adjusted
    PLR statistic of Theorem~\ref{thm:likl_ratio_2} for a nested model
    which excludes the first ten covariates for which
    $\bbeta_{0, j} = 0$, and $\chi^2_{10, q}$ is the $100 q \%$ quantile
    of a $\chi^2_{10}$ distribution. The constants $\bs$ and $\sigmas$
    result from the solution of~(\ref{eq:state_evol_inter}) for
    $\gamma^2 = 5$, and each combination of
    $\kappa \in \{0.2, 0.4, 0.6, 0.8\}$, and
    $\theta_0 \in \{0.5, 1.5, 2.5\}$.}
  \label{tab:intercept_LLR}
  \begin{center}
    \begin{small}
    \begin{tabular*}{\textwidth}{l@{\extracolsep{\fill}}l@{\extracolsep{\fill}}*{9}{S[table-format=-2.3]@{\,}}}
      \toprule
      & & \multicolumn{9}{c}{$100q$} \\
      \cmidrule(lr){3-11} 
      $\kappa$ & $ \theta_0 $ &  1  &  5  &  10  &  25  &  50  &  75  &  90  &  95  &  99 \\
      \midrule
      & 0.5 &  1.0 & 5.3 & 10.2 & 24.3 & 48.6 & 75.0 & 90.2 & 94.7 & 98.9 \\
       0.2  &  1.5  &  1.1 & 5.2 &  9.6 & 24.4 & 50.6 & 74.5 & 89.5 & 94.5 & 98.7 \\
      &  2.5  & 1.1 & 5.0 & 10.2 & 25.1 & 50.6 & 75.1 & 90.2 & 94.9 & 99.0 \\
      \midrule 
      &  0.5  & 1.0 & 4.9 &  9.5 & 24.7 & 50.8 & 75.2 & 89.6 & 95.0 & 98.8 \\
       0.4  &  1.5  & 1.0 & 5.4 & 10.7 & 25.6 & 49.4 & 75.0 & 90.0 & 94.6 & 99.0   \\
      &  2.5  & 1.3 & 5.3 &  9.9 & 25.3 & 49.4 & 75.4 & 90.4 & 95.2 & 98.9  \\
      \midrule 
      &  0.5  & 1.0 & 4.9 &  9.2 & 24.0 & 49.5 & 74.1 & 89.8 & 94.9 & 99.1  \\
       0.6  &  1.5  & 1.3 & 5.2 & 10.1 & 25.1 & 49.2 & 74.1 & 89.6 & 94.7 & 99.0   \\
      &  2.5  & 0.9 & 5.1 & 10.0 & 25.4 & 49.6 & 75.3 & 90.2 & 95.1 & 99.2  \\
      \midrule 
      &  0.5  & 1.4 & 5.3 & 10.1 & 24.5 & 50.3 & 74.5 & 89.5 & 94.8 & 99.0  \\
       0.8  &  1.5  & 0.7 & 4.9 &  9.3 & 24.2 & 49.8 & 74.7 & 90.3 & 94.9 & 98.9   \\
      &  2.5  & 0.9 & 4.7 &  9.6 & 24.8 & 49.8 & 75.1 & 90.4 & 95.0 & 99.1  \\
      \bottomrule
    \end{tabular*}
  \end{small}
\end{center}
\end{table}

To provide empirical support for the conjecture, we fix $n = 2000$ and
let $p$ vary so that $\kappa \in \{0.2, 0.4, 0.6, 0.8\}$, set
$\alpha = 1 / (1 + \kappa)$, and consider intercept values
$\theta_0 \in \{0.5, 1.5, 2.5\}$. We let $\bX$ be the $n$ times $p$
matrix with i.i.d $\mathrm{N}(0, 1/n)$ entries. Hence,
$\tau_j = 1/\sqrt{n}$. We also assume that $\bbeta_{0}$ has all but
$p / 4$ randomly chosen entries set to zero. The nonzero entries are
all of the same value and scaled such that
$\gamma^2 = \vnorm{\bbeta_{0}}_2^2 / n = 5$. For each combination of
$(\kappa, \theta_0)$, we draw $5000$ independent copies of
$\{\by,\bX\}$ from a logistic regression model with intercept. We
compute the adjusted PLR statistic from a nested model which
excludes the first ten covariates for which the corresponding signal
is zero. Table~\ref{tab:intercept_LLR} shows the estimated probability
that the adjusted PLR test statistic is less than the $100q\%$
quantile of a $\chi^2_{10}$ distribution, with $q$ ranging from $0.01$
to $0.99$. Clearly, the empirical distribution of the adjusted PLR
statistic closely agrees with the reference $\chi^2$ distribution.

In Section~\ref{subsec:adj_z_intercept} of the Supplementary Material
document, we contrast the adjusted $Z$-statistics for the first zero
(Table~\ref{tab:intercept_Zstat}) and first non-zero
(Table~\ref{tab:intercept_Zstat_nonnull}) coordinate of $\bbeta_0$ to
the $\mathrm{N}(0,1)$ quantiles. The empirical distributions of the
$Z$-statistics closely agree with the standard normal distribution.

\section{Estimation of unknown constants}
\label{sec:unknown_constants} 

\subsection{Without intercept} 

In the logistic regression model without intercept, the constants
$\kappa, \gamma$, which enter the system of equations
\eqref{eq:state_evol}, are unknown, and must be estimated.  A
straightforward estimate of $\kappa$ is $p / n$. For the estimation of
$\gamma$, we adapt the Signal Strength Leave-One-Out Estimator (SLOE)
estimator of \citet{yadlowsky+etal:2021}, which can be used to
estimate the corrupted signal strength
$ \upsilon^2 = \lim\limits_{n \to \infty}\var(\bx_i^\top \hat{\bbeta})
= \mu^2 \gamma^2 + \kappa \sigma^2$, whenever the ML estimator
$\hat\bbeta$ asymptotically exists.  Re-expressing the matrix
expressions in \citet{yadlowsky+etal:2021} into more familiar
statistical quantities that are readily available in the output of ML
estimation routines, the SLOE estimator of $\upsilon$ is defined as
\begin{equation}
  \label{eq:sloe}
  \hat\upsilon^2 = \frac{\sum_{i = 1}^n(s_i - \bar{s})^2}{n} \quad \text{with} \quad
  s_i =  \hat\eta_i - \frac{h_i}{1 - h_i} \frac{y_i - \zeta'(\hat\eta_i)}{ \zeta'(\hat\eta_i) (1 -  \zeta'(\hat\eta_i))} \,,
\end{equation}
where $\bar{s} = \sum_{i = 1}^n s_i / n$,
$\hat\eta = \bx_i^\top \hat\bbeta$, and $h_i$ is the $j$th diagonal
element of the ``hat'' matrix
$\bX(\bX^\top\bW(\hat\bbeta)\bX)^{-1}\bX^\top
\bW(\hat\bbeta)$. \citet{yadlowsky+etal:2021} show that
$ \hat{\upsilon}^2$ converges to $\upsilon^2$ in probability as
$n \to \infty$ in the setting of
Theorem~\ref{thm:aggregate_arbitrary}. Their proof relies on the ML
version of Theorem~\ref{thm:aggregate_arbitrary} and leave-two-out
techniques introduced by \citet{el_karoui:2018} and
\citet{sur+candes:2019}. We expect the same arguments to also work
when we replace $\hat\bbeta$ with $\betady$ in~(\ref{eq:sloe}). To get
estimates of $\mus, \bs, \sigmas$, we reparameterise the system of
equations \eqref{eq:state_evol} in terms of $\upsilon$ using the
identity $\gamma^2 = (\upsilon^2 - \kappa \sigma^2) / \mu^2$. Then,
replacing $\upsilon$ and $\kappa$ by their estimates, we can solve
\eqref{eq:state_evol} to obtain estimates for $\mus, \bs,
\sigmas$. The estimation of $\tau_j$ in
Theorem~\ref{thm:aggregate_arbitrary} and Theorem~\ref{thm:z_stats}
can be done using the residual sums of squares from the regressions of
each covariate to all the others, as detailed in
\citet[Section~5.1]{zhao+etal:2022}. The suitability of the above
estimates is supported by strong empirical evidence such as those
Section~\ref{subsec:corrected_ls} about the performance of 
non-oracle versions of the rescaled MDYPL estimator and of associated pivots.

\subsection{With intercept} 

When an intercept is included in the model, the set of unknown
parameters is $\kappa, \gamma, \theta_0$. We can estimate $\kappa$ as
$p / n$. \citet[][footnote~4, page 3]{yadlowsky+etal:2021} observe
that the estimator of $\upsilon^2$ in~(\ref{eq:sloe}) still converges
to $\upsilon^2 = \mu^2 \gamma^2 + \kappa \sigma^2$ even when an
intercept is included. We expect that the same holds when
$\hat{\bbeta}$ is replaced by $\betady$ in~(\ref{eq:sloe}).  Hence, we
can reparameterise \eqref{eq:state_evol_inter} in terms of $\upsilon$
instead of $\gamma$ and use the estimate of $\upsilon$ as we did
without intercept. The remaining unknown constant is
$\theta_0$. Rather than finding an estimate of $\theta_0$, we propose
to treat $\iota$ as a fixed parameter and solve
\eqref{eq:state_evol_inter} for $\mu, b, \sigma, \theta_0$ given
$\kappa, \eta, \iota$. According to Conjecture~\ref{conj:intercept}
$\hat{\theta}_0^{\ttextrm{DY}} \overset{p}{\longrightarrow}
\iota_{*}$, and, so, we can use the $\hat{\theta}_0^{\ttextrm{DY}}$ as
an estimate for $\iota$. We have found that this approach is rather
accurate and more accurate than estimating $\theta_0$ directly, as is
proposed in \citet[][Section~7.2]{zhao+etal:2022}. The estimation of
$\tau_j$ can be done as in the case without intercept. The numerical
experiments in the case study of Section~\ref{sec:case-study} provide
strong support for those procedures for estimating the unknown
constants.

\section{Illustration: Describing ``7''}
\label{sec:case-study}

We illustrate the utility of the MDYPL estimator and the application
of the results in Theorem~\ref{thm:aggregate_asymptotics}, and
Theorem~\ref{thm:likl_ratio_2} through the analysis of the Multiple
Features data set from the UCI repository
\citep{multiple_features:1998}. The dataset consists of digits
($0$-$9$) extracted from a collection of maps from a Dutch public
utility. Two hundred $30 \times 48$ binary images per digit were
available, which have then been used to extract feature sets
\citep[see][for details, where that dataset is used for assessing the
performance of various classifiers for digit
recognition]{jain+etal:2000}. We focus on the contribution of two of
those feature sets in describing the digit ``7''. The feature sets are
$76$ Fourier coefficients of the character shapes, which are computed
to be rotation invariant \citep[see, for example,][Section
4]{van_breukelen+etal:1998}, and 64 Karhunen-Lo\`eve coefficients. So
the available data consists of $2000$ indicators, taking values 1 or
0, depending on whether the digit is ``7" or not, respectively, and
$2000$ corresponding vectors with the values of the $140$ features. We
randomly split the data into a training set of $1000$ observations and
keep the remaining $1000$ observations for testing.

Depending on the font, the level of noise introduced during the
digitisation of the utility maps, and the downscaling of the digits to
$30 \times 48$ binary images, difficulties may arise in discriminating
instances of the digit ``7'' to instances of the digits ``1'' and
``4''. In addition if only rotation invariant features, like the
available Fourier coefficients, are used, difficulties may arise in
discriminating instances of the digit ``7'' to instances of the digit
``2''. So, we would expect to find evidence against the hypothesis
that a logistic regression model with an intercept and only Fourier
coefficients is an as good description of the digit ``7'' as the model
with both Fourier and Karhunen-Lo\`eve coefficients.

\begin{figure}[t]
  \centering
  \begin{subfigure}{.6\linewidth}
    \centering
    \includegraphics[width = \linewidth]{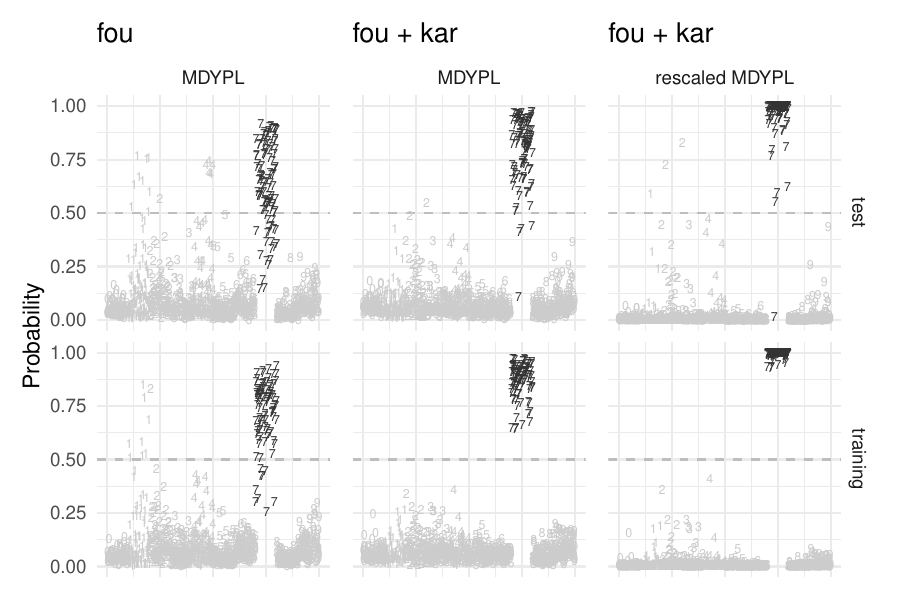}
    \caption{}
  \end{subfigure}%
  \begin{subfigure}{.39\linewidth}
    \centering
    \includegraphics[width = \linewidth]{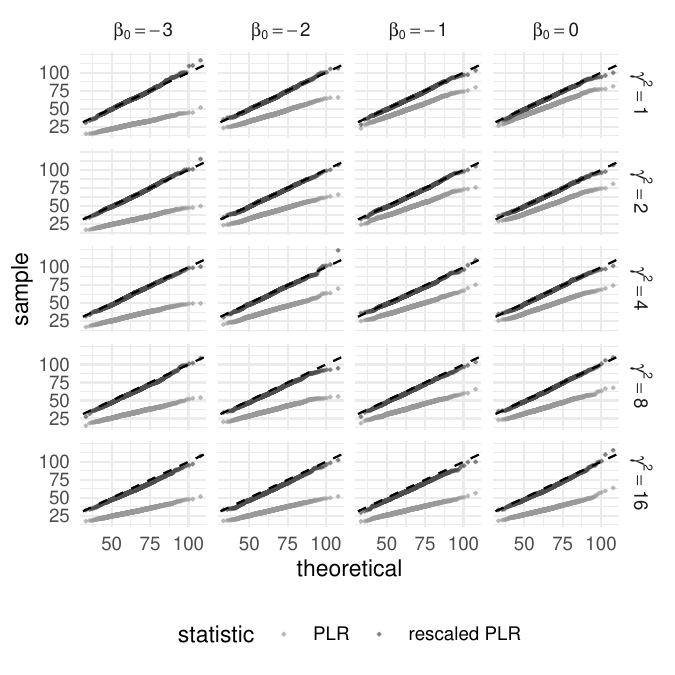}
    \caption{}
  \end{subfigure}%
  \caption{(a) Left: Estimated probabilities on the training and test
    sets, based on the MDYPL estimates of the model with Fourier
    coefficients only (\texttt{fou}). Right: Estimated probabilities
    on the training and test sets from the model with both Fourier and
    Karhunen-Lo\`eve coefficients (\texttt{fou+kar}), based on the
    MDYPL estimates (MDYPL), and their rescaled versions (rescaled
    MDYPL). The probabilities based on MDYPL are
    estimated as $1 / (1 + e^{-\etady_i})$ with
    $\etady_i = \theta_0^{\ttextrm{DY}} + \bx_i^\top \betady$, where $\bx_i$ are the
    Fourier and/or Karhunen-Lo\`eve coefficients for the $i$th digit
    in the training or test data, and $\betady$ are the MDYPL
    estimates of their coefficients. The probabilities based on the
    rescaled MDYPL estimator are estimated as
    $1 / (1 + e^{-\setady_i})$, where
    $\setady_i = \tilde{\beta}_0 + \bx_i^\top \betady / m$, where $m$ and
    $\tilde{\beta}_0$ are the estimates of $\mus$ and the intercept
    parameter from the procedure of
    Section~\ref{sec:unknown_constants}. (b) Q-Q plots for the
    penalized likelihood ratio statistic and its rescaled version for
    comparing models \texttt{fou+kar} and \texttt{fou}, using a
    $\chi^2_{64}$ reference distribution. }
  \label{fig:mfeat-fit}
\end{figure}

We test that hypothesis using the likelihood ratio test statistic, the
PLR statistic using $\alpha = 1 / (1 + p / n) = 1 / 1.14$, and the
rescaled version of it as suggested by
Theorem~\ref{thm:likl_ratio_2}. The constants $\bs$ and $\sigmas$,
required for the rescaled PLR statistic are estimated using the
procedure outlined in Section~\ref{sec:unknown_constants} using the
training data. Both the model with Fourier coefficients, and the model
with both coefficient sets, resulted in infinite ML estimates, and
perfect fits on the training data. As a result, the standard
likelihood ratio statistic has value $0$ providing no evidence
whatsoever against the hypothesis that the model with only Fourier
coefficients is adequate. Moreover, due to infinite estimates, the
methods of \citet{zhao+etal:2022} cannot be used to examine relevant
dimensionality corrections to the statistic. The PLR statistic takes
value $64.36$, which, when compared to the quantiles of a $\chi^2$
distribution with $64$ degrees of freedom, provides again no evidence
against the hypothesis. In stark contrast to those conclusions, the
rescaled PLR statistic takes value $173.34$,
and results in strong evidence in favour of the model with both
coefficient sets.

Figure~\ref{fig:mfeat-fit} (a) shows the estimated probabilities on
the training and test sets, based on the MDYPL estimates of the model
with Fourier coefficients only (\texttt{fou}), and the estimated
probabilities on the training and test sets from the model with both
Fourier and Karhunen-Lo\`eve coefficients (\texttt{fou+kar}), based
on the MDYPL estimates and their rescaled versions. As expected, model
\texttt{fou} has difficulties in discriminating between instances of
1, 2, 4, 7 in both the training and test data. The model \texttt{fou+kar} does better with MDYPL estimates, and markedly better when
based on $\tilde\theta_0$ and $\betady / m$ on both training and test
data, where $m$ is the estimate of $\mus$ and $\tilde\theta_0$ is the
estimated intercept from the procedures of
Section~\ref{sec:unknown_constants}.

The scaling of the PLR statistic we used earlier has been developed
under the assumption of a random $\bX$ with normal rows, which is
unlikely to be satisfied for the model matrix of \texttt{fou+kar},
is based on Conjecture~\ref{conj:intercept} when an intercept
parameter is present in the model, and relies on estimating the
scaling parameters according to the procedures in
Section~\ref{sec:unknown_constants}. So, we examine the performance
the rescaled PLR statistic by simulating $1000$ simulated response
vectors for each combination of the intercept parameter
$\theta_0 \in \{-3, -2, -1, 0\}$ and
$\gamma^2 \in \{1, 2, 4, 8, 16\}$. For each combination, the
$1000 \times 140$ matrix $\bX_{\ttextrm{fou+kar}}$ of Fourier and
Karhunen-Lo\`eve coefficients in the training set is centred to have
zero column means and kept fixed across simulations. The corresponding
parameters
$\bbeta_{\ttextrm{fou+kar}} = (\bbeta_{\ttextrm{fou}}^\top,
\bbeta_{\ttextrm{kar}}^\top)^\top$ are specified by setting
$\bbeta_{\ttextrm{kar}}$ to a vector of zeros, sampling the elements
of $\bbeta_{\ttextrm{fou}}$ independently from a standard normal
distribution, and rescaling them appropriately to ensure that the
sample variance of
$\bX_{\ttextrm{fou+kar}} \bbeta_{\ttextrm{fou+kar}}$ is
$\gamma^2$. Figure~\ref{fig:mfeat-fit} (b) shows Q-Q plots of the
samples of the PLR statistic and its rescaled version; clearly the
distribution of the rescaled PLR statistic agrees closely with the
nominal $\chi^2_{64}$ distribution.

\section{Concluding remarks}
\label{sec:summary}

We have established a framework for the recovery of estimation and
inferential performance in high-dimensional logistic regression under
proportional asymptotics, with independent normal covariates of
arbitrary covariance structure, irrespective of the existence of the
ML estimate. We achieve this by characterising the asymptotic
behaviour of the MDYPL estimator and associated inferential
procedures. In contrast to corresponding findings about ML in
\citet{sur+candes:2019} and \citet{zhao+etal:2022}, our findings are
no longer confined to the narrow region of $(\kappa, \gamma)$ values
below the phase transition curve of \citet{candes+sur:2020}, where the
ML estimates asymptotically exist. We also show that our results
generalise those for ML, reducing to the latter for particular values
of the DY prior hyperparameters. We provide a conjecture, with strong
empirical support, on the validity of our results when the intercept
parameter is non-vanishing, and develop estimation methods for all
unknown constants involved in our procedures. The strong empirical
evidence we provide for estimation and inference with or without an
intercept parameter and with non-normal covariates is reassuring about
the usefulness of MDYPL across a wide range of applications.


On a technical note, Section~\ref{appendix:thm1} of the Supplementary
Material document provides a rigorous approximation argument to apply
the AMP framework to the logistic regression model, which has been
missing from recent literature, as noted in
\citet{feng+et+al:2022}. That argument may be of independent interest
for the development of AMP instances where the Lipschitz assumption
fails. Finally, we establish asymptotic results for MPL with
non-separable penalty functions. This paves the way for similar
results for other penalty functions that imply a perturbation of the
responses, most notably the MJPL estimator.

\begin{figure}[t]
  \centering
  \includegraphics[width = \linewidth]{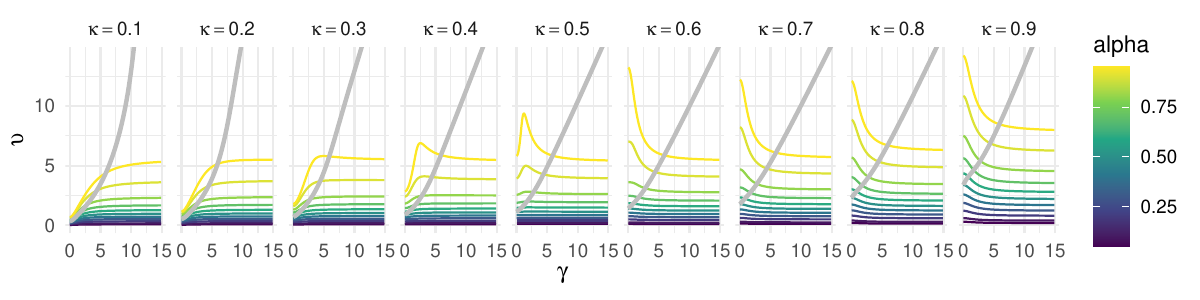}
  \caption{The function $\upsilon = h(\gamma)$ for every combination
    of $\kappa \in \{0.1, 0.2, \ldots, 0.9\}$ and
    $\alpha \in \{0.05, 0.1, 0.2, \ldots, 0.9, 0.95\}$ for
    $\gamma \in (0, 15)$. For each combination $\kappa$, $\alpha$, the
    functions are computed by solving the state evolution
    equations~(\ref{eq:state_evol}) without intercept over a dense
    grid of $\gamma$, and computing
    $\upsilon = \sqrt{\mus^2 \gamma^2 + \kappa \sigmas^2}$. The grey
    curve is the function $h(\cdot)$ for the ad-hoc choice
    $\alpha = 1 / (1 + \exp(-\gamma / 2))$.}
  \label{fig:upsilon-gamma}
\end{figure}

The estimation methods for the unknown constants rely on the
adaptation of the SLOE estimator $\hat{\upsilon}^2$ of
$\var(\bx_i^\top \betady) \to \upsilon^2$ in~(\ref{eq:sloe}), which
may result in some issues with particular data sets. For one, for any 
fixed value of $\kappa \in (0,1)$, the function $\upsilon = h(\gamma)$ is bounded. 
Towards a heuristic explanation, note that for
any fixed $\alpha < 1$, $\var(\bx_i^\top \betady)$ is bounded from above by a constant that
does not depend on $\gamma$ because
$\vnorm{\betady}_2^2 / \sqrt{p}$ is almost surely bounded uniformly over $\gamma > 0$ (see Theorem~\ref{thm:supp_boundedness_betady}). 
As a result, if $\hat\upsilon$ ends
up outside the image of $h(\cdot)$, it is not possible to find a
corresponding value for $\gamma$ to be used in the reparameterisation
of systems~(\ref{eq:state_evol}) and~(\ref{eq:state_evol_inter}) for
the procedures in Section~\ref{sec:unknown_constants}. 
Also, if
$\hat\upsilon$ lies in the image of $h(\cdot)$, there may be more than
one value of $\gamma$ that corresponds to $\hat\upsilon$. From the
consistency of $\hat\upsilon^2$ in~(\ref{eq:sloe}), we expect that the
former issue may only be present in small samples. For example,
Figure~\ref{fig:supp-upsilon-cls-vs-MDYPL} in the Supplementary
Material document shows the distribution of $\hat\upsilon$ from the
experiments of Section~\ref{subsec:corrected_ls}. Although we did not
encounter any issues with $\hat\upsilon$ outside the image of
$h(\cdot)$ when carrying out that simulation experiment, we see that
this may be possible if we took larger simulation sizes for $n = 400$
and that the probability of being outside the image decreases fast
with $n$. Figure~\ref{fig:upsilon-gamma} shows the function $h(\cdot)$
when there is no intercept in the model for a range of $\kappa$ and
$\alpha$. As is apparent, for any fixed $\alpha < 1$, $h(\cdot)$
transitions from an increasing function for small $\kappa$ to a
decreasing function for large $\kappa$, and we observe a range of
moderate $\kappa$, where two different $\gamma$ values result in the
same value for $\upsilon$. Clearly, the simulation study of
Section~\ref{subsec:corrected_ls} for $\kappa = 0.5$ results in an
$h(\cdot)$ function that is not invertible. However, all solutions to
the state evolution equations that we get have good properties, most
probably because we use the solution for the oracle choice of $\gamma$
as starting values. A solution to those issues, and the subject of
current work, may result from the use of an adaptive $\alpha$ in the
state evolution equations, where $\alpha$ is a function of
$\gamma$. For example, Figure~\ref{fig:upsilon-gamma} shows $h(\cdot)$
for $\alpha = 1 / (1 + \exp(-\gamma / 2))$.

The simulation results in Section~\ref{sec:case-study} and the
empirical evidence in Section~\ref{sec:bt-mar} of the Supplementary
Material document on the estimation and inference from Bradley-Terry
models with competitions missing at random suggest that our results
apply beyond normal covariates. In that direction, the techniques of
\citet{bayati+etal:2015} and \citet{chen+etal:2021} can be useful to
extend to relax that assumption to sub-Gaussian covariates. Also,
\citet{li+sur:2024} demonstrates the breakdown of proportional
asymptotic results under Gaussian model matrices for ridge linear
regression when the distribution of the covariates is heavy-tailed or
asymmetric, and develop appropriate debiasing procedures. The
extension of such techniques to MDYPL estimation with logistic
regression models is the subject of ongoing work. We have also assumed
the existence of a solution to equations \eqref{eq:state_evol} for
specific values of $(\alpha, \kappa, \gamma)$. Recent works such as
\citet{sur+candes:2019} and \citet{salehi+et+al:2019} also rely on
similar assumptions. The existence of stationary points for the state
evolution equations is proven in various works, for example,
\citet{montanari+etal:2023}, \citet{li+sur:2024}, and
\citet{bellec+loriyama:2024}, in settings that are similar in spirit
to ours. We expect the techniques used therein to be useful in
establishing the existence of a solution to \eqref{eq:state_evol} and
\eqref{eq:state_evol_inter} and its extensions to adaptive choices of
$\alpha$. Another line of ongoing work focuses on extending the
results about the MDYPL estimator to the case
$\bbeta_p \neq \bb{0}_p$, and proving Conjecture~\ref{conj:intercept},
which can be achieved using Gaussian min-max comparisons similar to
the analyses of \citet{salehi+et+al:2019}.

\section{Supplementary Materials}

The repository \url{https://github.com/psterzinger/MDYPL} provides the
Supplementary Material document and scripts to reproduce all analyses
and outputs in that document and the main text.


\section{Declarations}

We thank Federico Boiocchi for his careful reading of the manuscript
and useful suggestions during his internship at University of Warwick
in summer 2025, and for his substantial contribution to the
implementation of all the methods we developed in the brglm2
\citep{brglm2} R package.

For the purpose of open access, the authors have applied a Creative
Commons Attribution (CC BY) licence to any Author Accepted Manuscript
version arising from this submission.

\bibliographystyle{jss2}
\bibliography{mDYPL}

\includepdf[page=-]{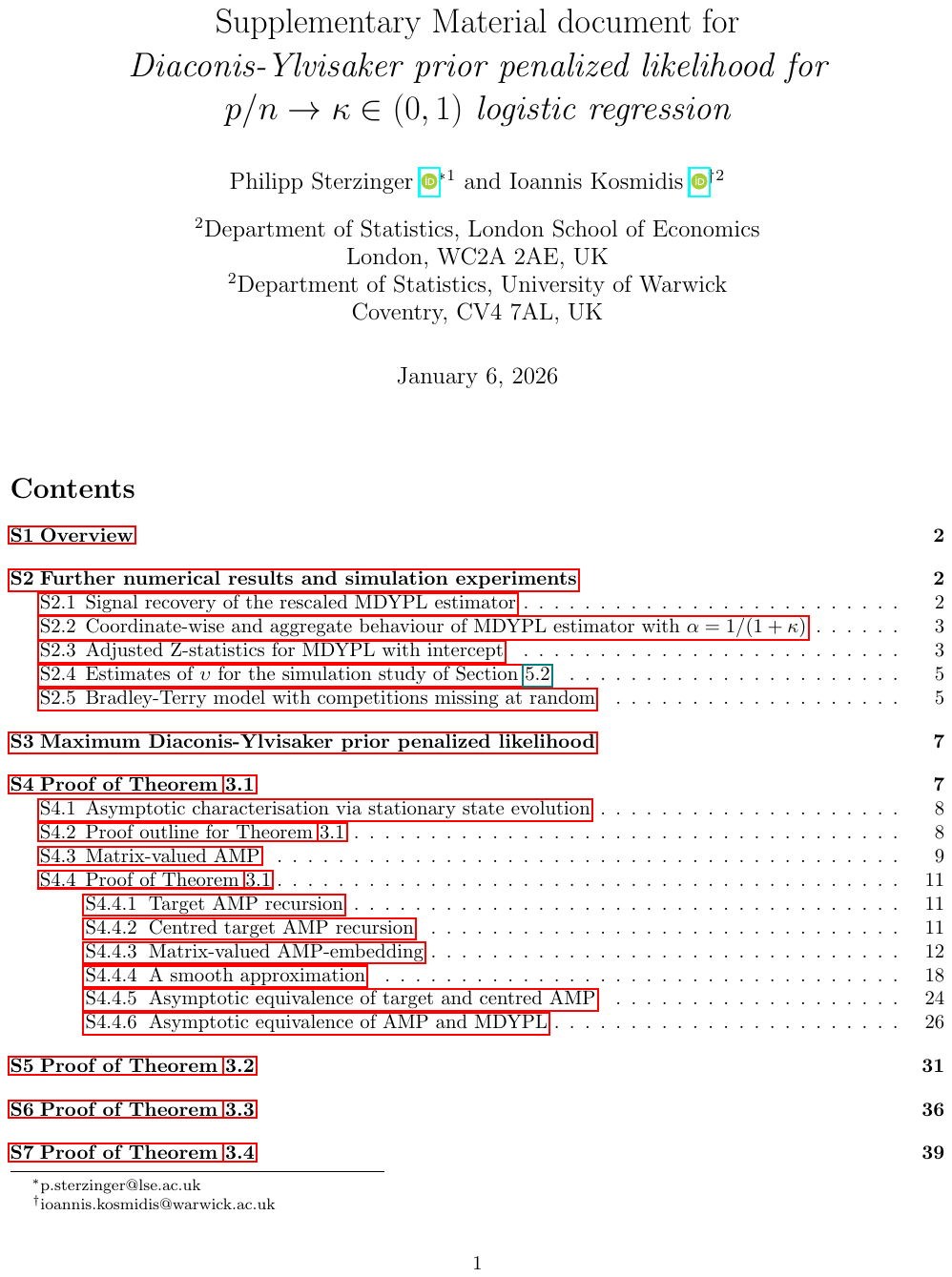}

\end{document}